\documentclass[11pt]{article}

\usepackage{ifpdf}
\usepackage{amsmath}
\usepackage{makeidx}
\usepackage{geometry}                
\usepackage{graphicx}
\usepackage{amssymb}
\usepackage{epstopdf}
\usepackage{lscape}
\usepackage{pdflscape}
\usepackage{verbatim}

\usepackage[gen]{eurosym}

\newcommand{\pf}{\vs \noindent {\it Proof:} \quad}

\newcommand{\C}{\mbox{\msbm{C}}}

\newcommand{\Q}{\mbox{\msbm{Q}}}
\newcommand{\R}{\mbox{\msbm{R}}}
\newcommand{\Z}{\mbox{\msbm{Z}}}

\renewcommand{\Re} {{\rm Re\,}}

\newcommand{\sinc}{\; {\rm sinc}\,}

\newcommand{\be}{\begin{equation}}
\newcommand{\ee}{\end{equation}}
\newcommand{\ba}{\begin{array}}
\newcommand{\ea}{\end{array}}
\newcommand{\bea}{\begin{eqnarray}}
\newcommand{\eea}{\end{eqnarray}}
\newcommand{\bean}{\begin{eqnarray*}}
\newcommand{\eean}{\end{eqnarray*}}

\newcommand{\vs}{\vspace{.25in}}

\newcommand{\bc}{\begin{center}}
\newcommand{\ec}{\end{center}}
\newcommand{\aster}{\begin{center} *** \end{center}}

\font\msbm=msbm10

\def\picture #1 by #2 (#3){
  \vbox to #2{
    \hrule width #1 height 0pt depth 0pt
    \vfill
    \special{picture #3} 
    }
  }

\def\scaledpicture #1 by #2 (#3 scaled #4){{
  \dimen0=#1 \dimen1=#2
  \divide\dimen0 by 1000 \multiply\dimen0 by #4
  \divide\dimen1 by 1000 \multiply\dimen1 by #4
  \picture \dimen0 by \dimen1 (#3 scaled #4)}
  }

\begin{document}

\title{
Function Theory on Inflationary Tessellations \\
\vspace{16pt}
\author{
{\sc Howard L. Resnikoff}\footnote{Resnikoff Innovations LLC; howard@resnikoff.com}
}   
\date{}
}

\maketitle

\begin{center}
In memory of \\
{\bf Max Koecher}
\end{center}

\vspace{12pt}

\begin{abstract}
This paper studies function theory on periodic and aperiodic inflationary tessellations. 
 
\end{abstract}


\newpage

\tableofcontents

\newpage

\section{Introduction}
\label{sec: introduction}

C. G. J. Jacobi called December 23, 1751 the birthday of the theory of elliptic functions \cite{Ayoub 1984}, but birth is preceded by conception. In this case, conception occurred 33 years earlier, in 1718, when Giulio Carlo Fagnano published three short papers concerned with measuring the length of an arc of the lemniscate  \cite{Fagnano}. Among other things, Fagnano discovered how to double the arc: 
\begin{small}
\begin{quote}
``Given the two equations
\bea
\frac{u \sqrt{2} }{\sqrt{1-u^4} }  &=& \frac{1}{z} \sqrt{ 1- \sqrt{1- z^4} } \\
\frac{dz}{\sqrt{1- z^4} } &=& \frac{2 \, du}{\sqrt{1-u^4} },
\eea
I say that if the first holds, so does the second."
\footnote{Fagnano's words, rendered by Ayoub \cite{Ayoub 1984}.}
\end{quote}
\end{small}

Carl Ludwig Siegel adopted Fagnano's work as the starting point for his three volume introduction to topics in complex analysis \cite{Siegel C L 1969}. He proposed a reconstruction of how Fagnano might have been led to his proof, and also simplified the argument, making it more elegant, and transforming the formulae into a more modern form, which is 
\be
 r^2 = \frac{2\, i\, u^2}{1-u^4}, \qquad \int_0^r \frac{dr}{\sqrt{1-r^4} } = (1+i)\, \int_0^u \frac{du}{\sqrt{1-u^4} }.
 \label{Fagnano1}
 \ee
Here we recognize an example of  the complex multiplication of elliptic curves.  David Hilbert said that the theory of complex multiplication of elliptic curves is not only the most beautiful part of mathematics, but of all science.\footnote{In 1932, at the International Congress of Mathematicians. \cite{Reid Constance}, p.200.} No wonder Siegel chose to begin with it.
 
The substitutions $r=\sqrt{2t}$,  $u= \sqrt{2 s} $ convert  eq(\ref{Fagnano1}) to
 \be
t= \frac{2  s}{\sqrt{4 s^2 - 1} }, \qquad
	 \int_0^{\sqrt{2t}} \frac{dt}{\sqrt{4t^3-t} } = (1+i)  \int_0^{\sqrt{2 s}} \frac{ds}{\sqrt{4s^3-s} } .
\label{lemn p-fncn cm}
 \ee
This leads us to the Weierstra{\ss} $\wp$-function as the uniformizing variable, which satisfies the differential equation $\wp'^2 =  4 \wp^3 - g_2 \, \wp - g_3 \,  \wp$; here,  $g_2=1$, $g_3=0$.  Hence the Klein modular invariant  $j(\tau)$ has the value $ j(1+i)= j(i)=1$, so the lattice of periodicity $\Omega$ for $\wp(z)$ is generated by 1 and $i$. Inverting the integrals, we find an expression for $\wp \left((1+i) z \right)$; this is a complex multiplication of the Weierstrass $\wp$-function. Amongst other things, it implies that $ (1+i) \Omega \subset \Omega$.  Equation (\ref{lemn p-fncn cm}) yields the complex multiplication  formula\footnote{We shall recover this formula   in a different way below. Cp. eq(\ref{-1+i complex mult}) et seq.}
\be
\wp \left((1+i) z \right) = -\frac{i}{8} \left( \frac{\wp'(z) }{\wp(z) } \right)^2.
\label{Fagnano complex mult}
\ee

\aster 

Every elliptic function is doubly periodic with respect to a lattice in $\C$  of rank 2, and every such lattice $\Omega$ admits multiplication by a rational integer. That is, for $n \in \Z^*$, $n \Omega \subset \Omega$.  
Lattice-invariant functions express this property through a real multiplication formula. For instance, for any lattice, the duplication formula for $\wp(z)$ corresponding to the inclusion $2\, \Omega \subset \Omega$  is
 $$ 
 \wp(2z) = \frac{1}{4} \left( \frac{ \wp''(z) }{  \wp'(x) } \right)^2 - 2 \wp(z) ,
 $$
which, after $\wp''$ and $\wp'$ are expressed  in terms of $\wp$ by means of the differential equation, yields a  formula that expresses  $\wp(2z)$ as a rational function of $\wp(z)$:
\be
\wp(2z) = \frac{1}{16 \wp(z)}  \frac{(12 \wp(z)^2 - g_2 )^2}{(4 \wp(z)^2 - g_2 ) }  - 2 \wp(z).
\label{wp-dup simplified}
\ee

This situation  is reminiscent of the duplication formula for $\sin z$, which can be expressed in four ways, viz, 
\bean 
\sin 2 z &=& 2 \sin z \cos z  \\
 &=& 2 \sin z \, \frac{ d \sin (z) }{dz}\\
   &=&   2 \sin z \, \sqrt{1 - \sin^2 z}\\
     &=&  2 \sin z \,  \sin (z + \pi/2) .
\eean
Although every non-zero rational integer is a multiplier,  it is unusual for a lattice to admit a complex multiplier, a complex number $\mu$ such that $\mu \,  \Omega \subset \Omega$. The necessary and sufficient condition is that $\Omega$ be the ring of integers of an imaginary quadratic  algebraic number field, and then $\mu$ can be any non-zero  integer in the field. In the lemniscatic case, $\Omega$ is the ring of Gaussian integers $\Q[i]$  in the field $\Q(i)$. 

\aster
 
Superficially,  multiplication formulae seem to arise as limiting special cases of an addition formula, and the addition formula appears to emerge from the periodicity of the  functions with respect to an underlying lattice. From a more general perspective, the existence of a multiplication formula is a consequence of a  theorem of Siegel that asserts the algebraic dependence of $n+1$ meromorphic functions on an $n$-dimensional  compact complex manifold \cite{Siegel C L 1955}. Here, the manifold is the elliptic curve, that is, the complex torus, on which the meromorphic functions live, and $\C$ is its universal covering space. Given the algebraic dependence theorem, it follows that $f(z)$ and $f(\mu z)$ satisfy an algebraic equation for any meromorphic function $f$ defined on the manifold and any complex number $\mu$  that carries the lattice of periods induced on the covering space into itself. 

There is no a priori reason, however, to believe that multiplication formulae do {\it not} exist in situations where meromorphic functions do not live on a compact complex manifold. Of course, some structure must be presupposed in order to prove anything. We have in mind inflationary tessellations -- particularly, aperiodic tessellations of the kind exemplified by Penrose tilings -- as the geometric substratum that substitutes for the lattice of periods in the classical theory, and functions defined in a way that mimics definitions of some of the most important classical higher transcendentals. 

One purpose of this paper is to construct multiplication formulae using the properties of positional representation for numbers. When the radix\footnote{Often called the `base' of the number system.} is a positive integer, the results recapitulate multiplication formulae for many classical functions. If the radix is not a rational integer, and particularly if it is a complex number that is not a quadratic imaginary integer, both the functions and the multiplication formulae appear to be new. 

Although these results may be interesting, they cannot be said of themselves to be important. Some of the new functions are similar to functions that appear in the theory of elliptic functions and higher arithmetic, but analogues of the properties that make those functions  significant are thus far unknown.  What may make the new functions, and  their multiplication formulae, of some immediate interest is their connection to inflationary tilings  and particularly to those that are aperiodic.

This observation sets the agenda for the paper.  We make use of the connection between inflationary tessellations and positional representation described in a previous paper \cite{HLR arXiv 3} to connect these tilings to positional representations, which can then be used to construct meromorphic functions on $\C$ that incorporate properties of the tiling structure. If, for instance, the inflationary tiling is a non-degenerate lattice in $\C$ (hence, periodic), then the functions are  elliptic functions.

Tessellations built up from an unlimited number of types of tiles are uninteresting.  As was shown in \cite{HLR arXiv 3}, if an inflationary tiling employs a finite number of types of tiles (the Penrose tiling uses  2), then  the inflation factor -- which turns out to be the radix of the associated positional representation -- must be an algebraic integer. Thus, although multiplication formulae exist for any radix, it is only the ones that are algebraic integers that couple function theory to an inflationary tessellation. This is another connection with, and an analog of,  complex multiplication. We shall begin by deriving various unrestricted  generalizations of classical multiplication formulae, but the heart of the paper concentrates on radices that are algebraic integers related to the golden number $\phi = \frac{1+\sqrt{5}}{2}$. In particular, we explore the radix $ i \sqrt{\phi}$. It  satisfies $x^4+x^2 = 1$ and is the multiplier for a particularly simple aperiodic tessellation of $\C$ that uses four types of tiles, three of which are geometrically similar. We call this the {\it special tessellation}.

\aster

The paper is organized to develop in detail this specific example of an aperiodic tiling of the plane and a family of associated functions that have complex multiplication formulae of the new kind.  But first we examine some classical functions from the standpoint of positional representations for numbers in $\R$ and in $\C$. Section \ref{sec: positional representations}  introduces positional representations, both real and complex, in a form appropriate for extension to inflationary tilings. Section \ref{sec: general paradigm} introduces the positional representation approach for some higher transcendental functions, and applies it to derive generalizations of some classical duplication formulae. When the radix is 2, these generalizations specialize to the classical duplication formulae.  In particular, these ideas are applied to  1-dimensional real multiplication formulae such as  eq(\ref{cot rho dup}), for example, which generalizes the classical duplication formula for $\pi \cot \pi z$, a trigonometric function that is structurally a 1-dimensional version of the Weierstra{\ss} elliptic function $\wp(z)$. 

Positional representations fall into two main classes: those for which the `integers' of the representation relative to the radix $\rho$ have a unique representation, and those which do not.  In section \ref{sec: non-unique examples} we investigate two examples  for which representations of $\rho$-integers are not unique.    

One of the pleasant peculiarities of  our approach is that partition of integers according to the number of digits in a positional representation leads to exact formulae for some associated functions. In section \ref{numerical exact identities}, exact  identities for partial sums of certain meromorphic functions defined by infinite series are derived. Since the partial sums are finite, questions of convergence are (temporarily) avoided. The Riemann $\zeta$-function has a ready generalization to this setting. Section \ref{sec: zeta analytic continuation} considers the generalization when the radix is the golden number $\phi = \frac{1+\sqrt{5}}{2}$ and develops a formula that provides an analytic continuation to the complex plane. 

Next we make explicit the connection with aperiodic tilings. Section \ref{sec: cot0_-phi} applies these ideas to a `duplication' formula for an aperiodic (1-dimensional) tessellation of $\R$ associated with the radix $-\phi$. This is preliminary to the work of section \ref{sec: function theory}, which takes up the analogous problem for the special tessellation, an aperiodic  tiling of the plane associated with radix $i \sqrt{\phi}$.   

Section \ref{sec: convergence} examines the convergence of the previously introduced  infinite series.    The paper concludes with remarks collected in section \ref{sec:remarks}.

\aster

The principal specific results of the paper are  for the complex radix $\rho = i \sqrt{\phi}$, presented in section \ref{subsec: wp_(-phi)}:   the  exact duplication formula for partial sums of the $\rho$-Hurwitz $\zeta$-function in eq(\ref{finite sum p-fncn Hurwitz s dup}), and the corresponding identity for the  Weierstra{\ss} $\wp$-function  $\wp_{\rho}(z)$ in eq(\ref{rho-wp dup}). The poles of  $\wp_{\rho}(z)$  are the vertices of the special tessellation. 

Another formula of interest provides an analytic continuation for the $\rho$-Riemann $\zeta$-function, which specializes to the usual Riemann $\zeta$-function and appears to provide something new even in this case; cp. eq(\ref{analytic continuation: classical}).

\section{Positional representations}
\label{sec: positional representations}

Many classical examples of multiplication formulae can be derived from the following paradigm.  Given a fixed integer $ \rho >1$,  express an arbitrary non-negative integer $n$ in the form $k \rho +r$ with $r = n \bmod \rho$. Terms  that have the same remainder $\bmod \, \rho$ in a  sum or product representing a function are collected.  This, apart from questions of convergence, supplies the desired formula. 

The positional representation of  $n$  in radix $\rho$ also expresses these relationships. If the digits of the representation are the integers $\Delta = \{ 0, 1, \dots, \rho -1 \}$, then 
$$ n= \sum_{k \geq 0} \delta_k \rho^k, \quad \delta_k \in \Delta, $$
where the sum is finite. This procedure can generalized. The radix need not be positive, nor  an integer. Any $\rho \in \C$ such that $| \rho | >1$ can serve as a radix for a positional representation. The set of digits can be selected in many ways, but the number of digits is $\lceil \, | \rho|^{\dim( \rho)} \, \rceil$, where $\dim(\rho) = 1$ for  representation of real numbers and $\dim(\rho) = 2$  for representation of complex numbers.\footnote{Roughly speaking, $\dim(\rho)$ is the dimension of the remainder set $R$ defined by eq(\ref{remainder set R}) below.}   $\lceil  x \rceil$ is the least integer greater than or equal to  $ x$.

 In either event, for $z$ running through some set of positive measure (which may be all of $\R$ or $\C$), 
\be
z = \sum z_k \rho^{-k}, \qquad z_k \in \Delta,
\label{general PR}
\ee
where finitely many positive powers of $\rho$ appear in the sum. The series converges absolutely since it is majorized by a geometric series.

It will be technically convenient, and easier to express the ideas we want to emphasize, if we limit ourselves ourselves to radices for which the number of digits  is 2.  In this case, $\Delta$ has two elements, which we select to be 0 and 1. Then the digits $z_k$ in eq(\ref{general PR}) are bits. This case already captures the essentials, both of complications and of consequences, of an arbitrary radix. Since only two digits occur, the corresponding generalized multiplication formula will be referred to as a {\it $\rho$-duplication} formula.

The set
\be
R := \left\{ \sum_{k=1}^{\infty} \delta_k\, \rho^{-k} \right\}, \quad \delta_{k} \in \Delta, 
\label{remainder set R}
\ee
is naturally called the set of {\it $\rho$-remainders} for the positional representation. It is $R$, and certain natural subsets of $R$, that will appear as the representative tiles of the associated tessellation.

Let $\Delta[X]$, resp. $\Delta[\rho]$,  denote the set of polynomials in $X$, resp. $\rho$, with coefficients in $\Delta$. The elements of  $\Delta[\rho]$  are called {\it $\rho$-integers} because $\Delta[\rho]$ is the usual set of non-negative integers when $\rho$ is a positive integer greater than 1.   

 Let $\Delta_n[X]$, resp. $\Delta_n[\rho]$, be the subsets of polynomials in $X$, resp. polynomials evaluated at $\rho$, that have $n$ coefficients, the $n$-digit polynomials. The leading coefficient may be 0.   The multiplicative containment $X \Delta_n[X] \subset \Delta_{n+1}[X]$ is the essential ingredient for deriving multiplication formulae.
 
 $\Delta_n[X]$ has $2^n$ distinct elements but the number of  elements of $\Delta_n[\rho]$ that have distinct values may be less. The coincidences are the consequence of algebraic identities satisfied by $\rho$.  Suppose that $p(X)$ and $q(X)$ are relatively prime polynomials with coefficients from $\Delta$. If $p(\rho) = q(\rho)$, then $p$ and $q$ are distinct representations for the same number. In particular, we see that in order for coincidences to occur, $\rho$ is necessarily an algebraic integer.

Let us write $p \sim q$ if $p(\rho)=q(\rho)$. ``$\sim$" is an equivalence relation. The collection of equivalence classes $\Delta[X]/(\sim)$ is in 1-1 correspondence with the values assumed by the polynomials.  

Specific sets  $\Omega$ of representatives of $\Delta[X]/(\sim)$ play an important role in what follows.  The important point is that a  set of unique representatives -- not necessarily complete -- of  $\rho$-integers be available.  

It is no accident that we use a notation --  $\Omega$ -- that was previously used to denote a non-degenerate lattice in $\C$, because when the lattice has a complex multiplication, it is a  set of unique representatives of  integers of the imaginary quadratic number field associated with the complex multiplication.

It was shown in \cite{HLR arXiv 3} that the inflation factor -- the ``multiplier" -- for an inflationary tiling must be an algebraic integer, whence it follows that in the applications of positional representation to tessellations, some $\rho$-integers may have more than one positional representation; in other words, the equivalence relation will not be trivial.\footnote{We cannot claim that there will necessarily be coincidences because the minimal polynomial of $\rho$ may not have coefficients from $\pm \Delta$. Thus the representations of $\rho$-integers for $\rho = \sqrt{2}$ are unique.} In the realm of $\rho$-multiplication formulae for meromorphic functions, this has the consequence that some terms of infinite series summed over the set $\Delta[\rho]$ will be repeated sufficiently often to force divergence.  The obvious solution is to select a unique representative polynomial from each equivalence class, i.e., to sum over the elements of $\Omega$. This selection will depend on the properties of the radix. 

Given $\Omega$, let $\Omega^*$ denote the set of non-zero elements, and denote the subset of $n$-digit $\rho$-integers by $\Omega_n$. Elements of $ \Omega_n$ may have 0 as the leading digit.  It  will be helpful to have a notation for the $n$-digit $\rho$-integers whose leading digit is 1, namely
\be
\nabla_n := \Omega_n - \Omega_{n-1} , \quad n>1 .
\label{nabla_n}
\ee
The set $\Omega_1$ of 1-digit $\rho$-integers is exceptional; we put  $\Omega_1 = \nabla_1 =\{ 0, 1 \}$. The $\nabla_n$ are pairwise disjoint, and  $\Omega = \bigcup_{n=1}^{\infty} \nabla_n$.   Put  $\Omega^*_n = \Omega_n - \{ 0 \}$.

\section{The general paradigm}
\label{sec: general paradigm}

As an introduction to the method, we shall formulate the general approach and then specialize it to obtain examples for classical functions.  Afterwards we will head into the uncharted waters and  unknown functions associated with radices for which the representations of $\rho$-integers are not unique. Throughout, we assume  $\Delta = \{ 0, 1 \}$.

\subsection{The $\rho$-Hurwitz $\zeta$-function}

Assume that $\lceil | \rho |^{\dim(\rho)} \rceil =2$. Suppose that $p(\rho) \neq q(\rho)$ for distinct polynomials $p,q  \in \Delta[X]$, so that $\Omega = \Delta[\rho]$.

The condition of uniqueness -- i.e., that each $\rho$-integer is represented exactly once in $\Delta[\rho]$ -- is satisfied if $\rho$ is a positive integer. It is also satisfied if $\rho$ is transcendental, for then the only solution to $p(\rho) = q(\rho)$  is  $p=q$.  

Algebraic cases sometimes lead to uniqueness.  Here are two examples: 

(1) If $\rho =3/2$, then  the only solution to $p(\rho) =  q(\rho) $ is  $p=q$. Suppose otherwise. After eliminating common powers of $\rho$  and clearing the denominator, the leading term is a power of 3, and the trailing term a power of 2.

(2) If $\rho = \sqrt{2}$, then a $\rho$-integer can be written in the form $m + n\sqrt{2} $ with $m$ and $n$ non-negative integers. Again, $p=q$ is the only solution. 

Cases where the representation of integers is not unique are of greater interest because they may be related to tessellations. For this reason we shall explore $\rho = \phi = \frac{1+\sqrt{5}}{2}$, the `golden number', and $\rho = i \sqrt{\phi}$. Here $\phi^2=\phi + 1$. The ratio of unique $n$-digit $\phi$-integers to the number of polynomials in $\Delta_n[X]$ approaches zero as $n$ increases. Nevertheless, we can specify a subset of polynomials such that every $\phi$-integer has a unique expression, and let the sum in an infinite series run over that set.  This will be done in section \ref{phi-unique}.

With $\Omega$ denoting a set of unique representatives of $\Delta[ \rho ]/(\sim)$, let  $\Omega_n$ denote the subset of $\Omega$ of representatives of $\Delta_n[\rho]$.   Set $\Omega^* = \Omega - \{ 0 \}$ and $\Omega^*_n = \Omega_n - \{ 0 \}$.  Where the $\rho$-positional representation of integers is not unique, we will explicitly state how the elements of $\Omega$ are chosen. We will write $\#(S)$ for the number of  elements in a set $S$.

\aster

Define the  $\rho$-Hurwitz  $\zeta$-function, 
\be
\zeta_{\rho}(s,z) := \sum_{\omega \in \Omega } \frac{1}{(z+ \omega)^s} , 
\label{rho-Hurwitz def}
\ee
and the $\rho$-Riemann $\zeta$-function, 
\be
\zeta_{\rho}(s) := \sum_{\omega \in \Omega^*} \frac{1}{\omega^s}  .
\label{rho-Riemann def}
\ee
Discussion of convergence of the series will be postponed to section \ref{sec: convergence} because it depends on particular properties of the radix.

Suppose that the series converge absolutely and uniformly on compact subsets of the complement of $\Omega$, resp. $\Omega^*$, and that the representation of $\rho$-integers in $\Delta[ \rho ]$ is unique. Then 

\bean
\rho^s \, \zeta_{\rho}(s, \rho\, z) 
	&=& \sum_{\omega \in \Omega} \frac{1}{(z+ \omega/ \rho)^s} 			\nonumber \\
	&=& \sum_{\omega \in \Omega} \frac{1}{(z+ \omega)^s} 	 + 			\sum_{\omega \in \Omega} \frac{1}{(z+ \frac{1}{\rho} + 			\omega)^s} .  \nonumber \\
\eean
Hence

\be
\rho^s \, \zeta_{\rho}(s, \rho \, z)  =  \zeta_{\rho}(s, z) + \zeta_{\rho}(s,  z +1/\rho),
\label{gen Hurwitz dup}
\ee
which is the $\rho$-duplication formula for the $\rho$-Hurwitz $\zeta$-function.

The function $ \zeta_{\rho}(s,  z +1/\rho)$ is holomorphic in a neighborhood of $z=0$. From eq(\ref{gen Hurwitz dup}) we find
\bean
 \zeta_{\rho}(s,  1/\rho) 
 	&=& \lim_{z \rightarrow 0}  \left(  \rho^s \, \zeta_{\rho}(s, \rho\, z)  - \zeta_{\rho}(s,  z)  \right) \\
	&=& \lim_{z \rightarrow 0} \left(
		\left\{ \frac{1}{z^s} + \sum_{\Omega^*} \frac{ \rho^s}{(\rho z+ 		\omega)^s} \right\} -
		\left\{ \frac{1}{z^s} + \sum_{\Omega^*} \frac{ 1}{(\rho z+ 			\omega)^s} \right\} 
		\right)  \\
	&=& \left( \rho^s - 1 \right)\zeta_{\rho}(s).
 \eean
 Hence
 \be
\left( \rho^s - 1 \right)\zeta_{\rho}(s) =   \zeta_{\rho}(s,  1/\rho)  
  \label{gen Hurwitz>Riemann}
 \ee
expresses the $\rho$-Riemann $\zeta$-function in terms of the  $\rho$-Hurwitz $\zeta$-function.  We have implicitly used uniqueness to conclude that $-1/\rho$ is not a pole, for if it were, then it would have the same value as some $\omega \in \Omega^*$ and clearing the denominator would yield two representations for 1.

The poles of $\zeta_{\rho}(s,z)$ are the points  $-\Omega$.\footnote{It would be aesthetically more satisfactory to define the Hurwitz function, etc. so that the set of poles is $\Omega$, but,  for the moment, we prefer to make the connection with classical notations.}   We shall calculate the Laurent series for $\zeta_{\rho}(s,z)$ at $z=0$. In a neighborhood of $z=0$, the function
$$ h(z):= \zeta_{\rho}(z) - \frac{1}{z^s} =\sum_{\Omega^*} \frac{1}{(z + \omega)^s} $$
is holomorphic. Then
\bean
\left. d^k h/dz^k \right|_{z=0}
	&=& \left. (-1)^k \, \frac{\Gamma(s+k)}{\Gamma(s)} \, 				\sum_{\Omega^*} \frac{1}{(z + \omega)^{s+k}}  \right|_{z=0}\\
	&=&  (-1)^k \, \frac{\Gamma(s+k)}{\Gamma(s)} \, \zeta_{\rho}(s+k) .
\eean
Therefore
\be
\zeta_{\rho}(s,z ) = \frac{1}{z^s} + \sum_{k=0}^{\infty} (-1)^k \frac{\Gamma(s+k)}{\Gamma(s) \Gamma(k+1)} \, \zeta_{\rho}(s+k) z^k.
\label{zeta_rho-Hurwitz Laurent series}
\ee

Bearing mind that the constant term on the right side of  eq(\ref{zeta_rho-Hurwitz Laurent series}) is $ \zeta_{\rho}(s) $, evaluate both sides  at $z=1/\rho $ and use eq(\ref{gen Hurwitz>Riemann}) to find  
\be
\zeta_{\rho}(s) =  \frac{\rho^s}{\rho^s - 2} \,\left( 1+  \sum_{k \geq 1} (-1)^k \, \frac{\Gamma(s+k)}{\Gamma(s) \Gamma(k+1)} \, \frac{\zeta_{\rho}(s+k)}{\rho^{s+k}} \right).
\label{zeta_rho-Riemann Laurent series}
\ee
Equation (\ref{zeta_rho-Riemann Laurent series})  provides an analytic continuation of the $\rho$-Riemann $\zeta$-function to the entire complex plane.

If eq(\ref{zeta_rho-Riemann Laurent series}) is valid (i.e., for values where the series converges appropriately), the function has a simple pole at
\be
s = \log 2/\log \rho \, .
\label{rho-Riemann-zeta pole}
\ee

\subsection{Classical examples: $\rho \in \R$, $\Delta = \{ 0,1\} $}

Classical examples of duplication formulae fall into several classes. For a real radix, $ | \rho | = n \in \Z^+$: that is, an endomorphism of the lattice $\Z$.  The ancillary condition $\lceil | \rho | \rceil=2$ implies  $\rho = \pm 2 $. For $\rho \in \C$, the radix is an endomorphism of a non degenerate lattice in $\C$ and hence a rational or algebraic integer in a quadratic imaginary number field that depends on the lattice. In this case, the ancillary condition is  $| \rho |^2 =2$. There are just three essentially different cases, which we shall consider separately: $\rho = i \sqrt{2}$,  $\rho =\frac{1+i \sqrt{7}}{2}$, and $\rho = (-1 + i)$.\footnote{These multipliers are only determined up to a unit of the imaginary quadratic field.}  

We begin with $\rho = 2$.

\subsubsection{The Riemann $\zeta$-function}

When $\rho=2$, the $\rho$-Hurwitz and the $\rho$-Riemann functions specialize to the classical Hurwitz and Riemann $\zeta$-functions, respectively.  

Equation (\ref{gen Hurwitz>Riemann}) yields  
 $\zeta_{\rho}(s) = \zeta_{\rho}(s,1) = \zeta(s)$, where $\zeta(s)$ is the Riemann $\zeta$-function, which expresses $\sum_{k=0}^{\infty} \frac{1}{(2k+1)^s}$ in terms of $\zeta(s)$.  From eq(\ref{rho-Riemann-zeta pole}) we find the  pole of $\zeta(s)$ at $s=1$.

Equation(\ref{zeta_rho-Riemann Laurent series})  reduces to an identity for the Riemann $\zeta$-function, namely
\be
\zeta (s) =  \frac{2^s}{ 2^s - 2 } \left( 1 + \sum_{k=1}^{\infty} (-1)^k   \frac{\Gamma(s+k)}{\Gamma(s) \Gamma(k+1)}  \, \frac{\zeta (s+k)}{2^{s+k}}  \right) ,
\label{analytic continuation: classical}
\ee 
which is probably  known but I have not found it in the literature. If $s = \frac{1}{2} + it$ is a zero, it implies the curious property
\be
\sum_{k=1}^{\infty}  \frac{(-1)^k}{2^{k}}  \frac{\Gamma(k+1/2 + it )}{ \Gamma(k+1)}\, \zeta \left( k+1/2 + it \right)  = -2^{1/2 + it} \, \Gamma(1/2 + it) .
 \label{curious 1}
 \ee
 Moreover, if we temporarily write
\be 
g(s) := 1 + \sum_{k=1}^{\infty} (-1)^k   \frac{\Gamma(s+k)}{\Gamma(s) \Gamma(k+1)}  \, \frac{\zeta (s+k)}{2^{s+k}} ,
\ee
then, since the zeros of $2^s - 2 $ are $1 + \frac{2 \pi i n}{\log 2}, \, n \in \Z$, it follows that $g\left( 1 + \frac{2 \pi i n}{\log 2} \right) =0$. 

\aster 

\begin{figure}[h]
\begin{center}
\caption{Graph of the absolute value of the Riemann $\zeta$ function on the critical line using analytic continuation provided by eq(\ref{analytic continuation: classical}) -- the black curve -- compared with the standard function.}
\label{zeta graph}
\includegraphics[width=2.5in]{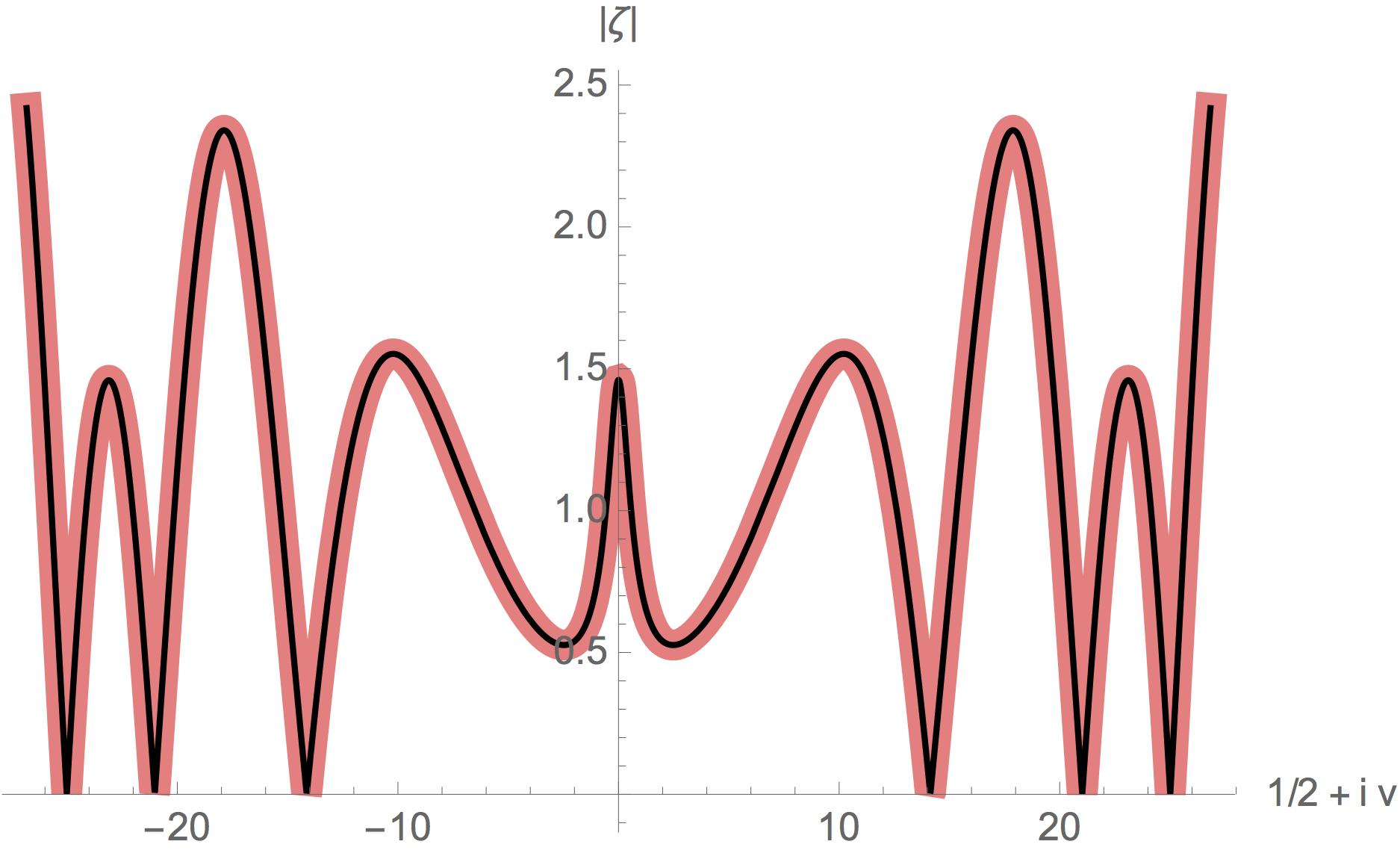}
\end{center}
\end{figure}

It may not be amiss to devote a few words to what  eq(\ref{analytic continuation: classical})  means in practice. The series converges for $\Re(s)>1$ and the $\zeta$-functions that appear in it are $\{ \zeta(s+k) : k \in \Z^+ \}$, so we can calculate  $\zeta(s)$ on the critical line $\Re(s)=1/2$ from the series. Figure \ref{zeta graph} shows the graph of $| \zeta(\frac{1}{2} + i v) |$ calculated two ways: the thick pink curve was computed by ${\it Mathematica}^{\copyright}$ using its internal code; the black curve was calculated using 80 terms of the series from eq(\ref{analytic continuation: classical}).

\subsubsection{The $\rho$-polygamma function}

The gamma function of Euler has the infinite product expansion
\be
\frac{1}{\Gamma(z)} = e^{\gamma z} \, z \, \prod_{k=1}^{\infty} \left( 1 + 		\frac{z}{k} \right) e^{-z/k} ,
\ee
where  $\gamma$  is Euler's constant.  Gauss' multiplication formula,  
\be
 \quad  \Gamma (n z) = (2 \pi)^{\frac{1-n}{2} } n^{n x -\frac{1}{2} } 
		\prod_{k=0}^{n-1}  \Gamma \left( z +\frac{k}{n} \right) , \quad n \in \Z^+,
\ee 
follows by splitting the infinite product into factors according to the remainders $\bmod \, n$. The duplication formula is the special case
\be
\Gamma(2 z) = \frac{2^{z-1} }{ \sqrt{\pi} }
		 \Gamma(z) \Gamma \left(z+\frac{1}{2} \right)  .
\ee  
We shall retrieve this duplication formula in the next subsection.

The polygamma functions, that is, the derivatives of $\log \Gamma(z)$, have the expansion
\bea
\psi^{(0)}(z)=\frac{d \log \Gamma(z) }{dz} &=& -\frac{1}{z} - \gamma -
	\sum_{k=0}^{\infty} \left( \frac{1}{z+k} -\frac{1}{k} \right), \\
\psi^{(n)}(z)=\frac{ d^{n+1}\log \Gamma(z) }{ dz^{n+1}} &=& (-1)^{n+1} n! \sum_{k=0}^{\infty} \frac{1}{(z+k)^{n+1} }, \quad 0<n \in \Z^+ .
\eea
Evidently $\psi^{(n)}(z)$  is proportional to the 2-Hurwitz $\zeta$-function, 
\be
\psi^{(n)}(z)= (-)^{n+1} n! \, \zeta_{2}(n+1,z), \quad n \in \Z^+,
\ee
so it too will inherit a $\rho$-duplication formula. 

Introduce the $\rho$-polygamma function $\psi_{\rho}^{(n)}(z)$ for $ n \in \Z^+$ by
\be
\psi_{\rho}^{(n)}(z) =  (-1)^{n+1} n! \, \zeta_{\rho}(n+1,z) .
\ee
Regarding $\psi_{\rho}^{(0)}(z)$, we must reason somewhat differently.  Introduce
\be
\zeta_{\rho}(1,z):= \frac{1}{z} + \sum_{\omega \in \Omega^*} \left( \frac{1}{z+ \omega} -  \frac{1}{ \omega} \right) .
\label{zeta_rho(1,z)}
\ee
Formally,
\bea
\rho \, \zeta_{\rho}(1,\rho z) 
	&=& \frac{1}{z} +  \sum_{\omega \in \Omega^*} \left(\frac{1}{z+\omega/\rho} 		- \frac{1}{\omega/\rho} \right) \nonumber  \\
 	&=& \frac{1}{z}  + \sum_{\omega \in \Omega^*} 								\left(\frac{1}{z+  \omega} - \frac{1}{  \omega} \right)  +  \sum_{\omega 		\in \Omega} \left( \frac{1}{z+1/\rho + \omega} - \frac{1}{1/\rho + 				\omega} \right)  \nonumber \\
	&=& \zeta_{\rho}(1, z)  + \frac{1}{z+1/\rho} -\frac{1}{1/\rho} +  \nonumber \\
	&&   \qquad \sum_{\omega \in \Omega^*} \left( \frac{1}{z+1/\rho + \omega} -		\frac{1}{\omega} +\frac{1}{\omega}   - \frac{1}{1/\rho + \omega}  \right)   		\nonumber \\
	&=& \zeta_{\rho}(1, z)  + 
		\frac{1}{z+1/\rho} + \sum_{\omega \in \Omega^*} \left( \frac{1}{z+1/\rho 		+  \omega}  -  \frac{1}{\omega}  \right) -  \nonumber \\
	&&   \hspace{58pt}   \frac{1}{1/\rho } - \sum_{\omega \in \Omega^*} \left( \frac{1}			{1/\rho +  \omega}  -  \frac{1}{\omega}  \right) \nonumber \\
\rho \, \zeta_{\rho}(1,\rho z) 
	&=& \zeta_{\rho}(1, z)  +  \zeta_{\rho}(1, z + 1/\rho)  - \zeta_{\rho}(1,  1/\rho) , 
\label{zeta_rho(1,z) dup}
\eea
which is the $\rho$-duplication formula for $ \zeta_{\rho}(1, z) $, and hence is proportional to that for $\psi_{\rho}^{(0)}(z)$.  The duplication formula for $\psi_{\rho}^{(n)}(z), \, n>0$ can also be derived by differentiating this formula.

Specializing to $\rho=2$, we find $\zeta_{2}(1/2) = 2 -  \sum_{k=1}^{\infty}\frac{1}{k(2k+1)}=2 \log 2$ so we recover the duplication formula for the classical  polygamma function:
\be 
2 \psi^{(0)}(2z) = 2 \log 2 + \psi^{(0)}(z) + \psi^{(0)}(z+1/2).
\ee

\subsubsection{The $\rho$-Gamma function}

Working backward, integration followed by exponentiation of eq(\ref{zeta_rho(1,z)}) leads to an infinite product we shall call the $\rho$-Gamma function:
\be
\frac{1}{ \Gamma_{\rho}(z) } := e^{\gamma_{\rho} z} \, z \, \prod_{\omega \in \Omega^*} \left( 1 + \frac{z}{\omega} \right) \, e^{-z/\omega} .
\label{rho=Gamma fncn}
\ee
The constant of integration $\gamma_{\rho}$ generalizes Euler's constant, to which it specializes for $\rho = 2$, in which case  $\Gamma_{2}(z) = \Gamma(z)$.

From the $\rho$-duplication formula eq(\ref{zeta_rho(1,z) dup}) we find the $\rho$-duplication formula for $\Gamma_{\rho}(z)$:
\be
\Gamma_{\rho}(\rho  z) = e^{- \psi_{\rho}^{(0)}(1/\rho) \, z } \, \Gamma_{\rho}(z) \,  \Gamma_{\rho} (z+ 1/\rho).
\ee 

\subsubsection{The $\rho$-cotangent function: $\rho = -2$}

The Laurent series for the cotangent function and its derivatives are
\bea
\pi \cot \pi z &=& \frac{1}{z} +
	 \sum_{k \in \Z-\{ 0 \} } \left( \frac{1}{z+k} - \frac{1}{k} \right) ,
 \label{cot0} \\
\frac{d^m \pi \cot \pi z}{dz^m} &=& (-)^n n!  \sum_{k \in \Z }  \frac{1}{(z+k)^{m+1}} , \qquad m>0.
\label{cotn}
\eea
The derivatives of $\pi \cot \pi z$ can be expressed in terms of the $\rho$-Hurwitz $\zeta$-function with $\rho=-2$.  Note that here the summation index set $\Omega$ of $(-2)$-integers  is $\Z$. To see this, write
\bean
\left\{  \sum_{k \geq 0} \delta_k (-2)^k , \quad \delta_k \in \Delta \right\} 
	&=& \left\{    \sum_{k \geq 0} \delta_{2k} 4^{k} -  \sum_{k \geq 0} (2 			\delta_{2k+1}) 	4^{k}  \right\}  \\
	&=& \left\{  \sum_{k \geq 0} \left( \delta_{2k}  - 2\, \delta_{2k+1} \right)4^{k} 		\right\}  \\
	&=& \left\{  \sum_{k \geq 0}   \epsilon_k  4^k \, : \, \epsilon_k \in  
	\{ -2,-1,0,1 \}   \right\}  \\
	&=& \Z ,
\eean
because   $\{ -2,-1,0, 1 \}$ is a maximal set of integers  inequivalent  modulo 4 and therefore suitable as a set of digits for a radix 4 positional representation. Comparison of eq(\ref{cotn}) with eq(\ref{rho-Hurwitz def}) shows that
\be
\frac{(-)^n}{n!} \frac{d^m \pi \cot \pi z}{dz^m} = \zeta_{-2}(m+1,z), \quad  m \in \Z^+
\ee

\aster

The duplication formula for $\pi \cot \pi z$  is obtained in the usual way, by integrating the series for its first derivative and selecting the constant of integration to force convergence. This also works for the $\rho$-cotangent function $\cot_{\rho} z$, which we define as follows, omitting the factors of $\pi$ for notational simplicity.\footnote{Thus $\cot_{-2}( z) = \pi \cot \pi z$.}

Suppose the series
\be
\cot_{\rho} (z) := \frac{1}{z} + \sum_{\omega \in \Omega^*} \left( \frac{1}{z + \omega} - \frac{1}{\omega} \right)
\ee
converges absolutely and uniformly on the complement of $\Omega$ for the given $\rho \in \R$ with $1< | \rho | \leq 2$.
Recalling that 
$$
\omega = \delta + \rho\,  \tilde{\omega}, \quad \delta \in \Delta; \quad  \omega, \tilde{\omega} \in \Delta[\rho] =  \Omega,$$
calculate
\bean
\rho \cot_{\rho} (\rho  z) 
	&=& \frac{1}{z} + \sum_{\omega \in \Omega^* }  			\left(   \frac{1}{z + \omega/\rho} - \frac{1}{\omega/\rho}  \right)  		\nonumber \\
	&=& \frac{1}{z} + 
		\sum_{\omega \in \Omega^*}  								\left(  \frac{1}{z  + \omega} - \frac{1}{  \omega}  \right)  + 			\sum_{\omega \in \Omega }  								\left(  \frac{1}{z+ 1/\rho  + \omega} - \frac{1}{ 1/\rho + \omega}  		\right)   \nonumber  \\
	&=& \cot_{\rho} ( z)   + 
		 \sum_{\omega \in \Omega}  									\left(  \frac{1}{z + 1/\rho + \omega} - \frac{1}{ 1/\rho + \omega}  \right)   \nonumber \\
	&=&  \cot_{\rho} ( z)   + 
		 \frac{1}{z + 1/\rho}  +
		 \sum_{\omega \in  \Omega^*  }  								\left(  \frac{1}{z + 1/\rho + \omega}    - \frac{1}{\omega}  \right) - \nonumber \\
	&&  \qquad \frac{1}{ 1/\rho} -
		 \sum_{\omega \in   \Omega^*  }  								\left(  \frac{1}{1/\rho + \omega}   - \frac{1}{\omega} \right) \nonumber \\
	&=&  \cot_{\rho}(z) +    \cot_{\rho}(z+1/\rho)  - \cot_{\rho} (1/\rho).
\label{rho-cot-repl}
\eean

Thus, with  conditions that insure  convergence, the $\rho$-duplication formula for $\cot_{\rho}(z)$ is
\be
\rho \cot_{\rho} (\rho  z)  = \cot_{\rho}(z) +   \cot_{\rho}(z+1/\rho)  - \cot_{\rho} (1/\rho).
\label{cot rho dup}
\ee

For $\rho =- 2$, $\cot_{-2}(z) =\pi \cot \pi z$ and eq(\ref{cot rho dup})  is the classical duplication formula.

\subsection{Classical examples: $\rho \in \C - \R$, $\Delta=\{0,1 \}$}

Now we turn to complex radices.  It has already been noted that there are three  lattices  that admit  complex multiplications $\rho$ for which $| \rho |^2=2$.  The associated quadratic imaginary field is $\Q(\rho)$. The lattice  generators can be taken as $\{ 1, \rho \}$ with $\rho  \in \{  i\sqrt{2},\frac{1 +i \sqrt{7} }{2},  -1+i \}$.  The radix $\rho = - 1 +i$  corresponds to the lemniscatic complex multiplication discussed in the introduction.\footnote{The unit $-1 \in \Z$ makes a difference in positional representations with real radix. For instance, $\Delta[2]= \Z^+ \cup \{ 0 \}$ whereas $\Delta[-2] = \Z$. An analogous statement is true for complex radices. The units for $1+i$ are $\{ \pm1, \pm i\}$.  Radix $1+i$ yields a subset of the Gaussian integers $\Q[i]$ which, together with $i- \Q[i]$,  exhausts all the Gaussian integers, whereas  radix $-1+i $ produces  every Gaussian integer: $\Delta[-1+i] = \Q[i]$. }  In each case the modular invariant is a rational number:
\be
j(i \sqrt{2}) = (5/3)^3, \quad j \left( \frac{1+i \sqrt{7}}{2} \right) = -(5/4)^3,   \quad  j(-1+i) = j(1+i) = j(i) =1.
\ee

\subsubsection{The $\rho$-Weierstra{\ss} $\wp$-function}

Recall the definition of the $\rho$-Hurwitz $\zeta$-function eq(\ref{rho-Hurwitz def}),
$$ \zeta_{\rho}(s,z):= \sum_{\omega \in \Omega} \frac{1}{(z+ \omega)^s}, \quad \Omega = \Delta[ \rho ], $$
and the $\rho$-duplication formula eq(\ref{gen Hurwitz dup}),
$$ \rho^s  \zeta_{\rho} (s, \rho z) =  \zeta_{\rho} (s,  z)  +  \zeta_{\rho} (s,  z+1/\rho) ,$$
which is valid for $\rho \in \C$ when the series converges.  Recall eq(\ref{gen Hurwitz>Riemann}),
$$
\zeta_{\rho}(s, 1/\rho) =    (\rho^s - 1) \zeta_{\rho}(s).
$$
Introduce the $\rho$-Weierstra{\ss} function
\be
\wp_{\rho}(z) := \frac{1}{z^2} + \sum_{\omega \in \Omega^*} \left( \frac{1}{(z+ \omega)^2} - \frac{1}{( \omega)^2}  \right) 
\label{rho-p-function}
\ee
and its derivatives,
\bea
\wp_{\rho}^{(k)}(z)  &:= & (-)^k \, (k+1)!  \, \sum_{\omega \in \Omega} \frac{1}{(z+ \omega)^{k+2} }  \nonumber \\
	&=&   (-)^k \, (k+1)! \, \zeta_{\rho}(k+2,z), \quad k \in \Z^+ .
\label{rho-p-function n deriv}
\eea

Assuming there are no duplications in $\Delta[\rho]$, repetition of the argument that led to eq(\ref{rho-cot-repl}) for the $\rho$-cotangent function yields the $\rho$-duplication formula for $\wp_{\rho}(z)$:
\be
\rho^2 \, \wp_{\rho}( \rho z) = \wp_{\rho}(  z) + \wp_{\rho}( z+1/\rho) - \wp_{\rho}( 1/\rho).
\label{rho wp dup}
 \ee

Regarding the Laurent expansion of $\wp_{\rho}(z) - 1/z^2 $ at $z=0$, , the constant term is 0 ($= \sum_{\Omega^*} \left( \omega^{-s}- \omega^{-s} \right)$ ) and eq(\ref{rho-p-function n deriv}) implies
$$ \left. \frac{d^k \left( \wp_{\rho}(z) - z^{-2}\right) }{dz^k} \right|_{z=0} = (-1)^k (k+1)! \, \zeta_{\rho}(k+2) $$
so
\be
\wp_{\rho}(z) = \frac{1}{z^2} + \sum_{k=0}^{\infty} (-)^k (k+1) \zeta_{\rho}(k+2)\, z^k  .
\label{rho-wp Laurent}
\ee
This formula shows that the function $\zeta_{\rho}(k+2)$ is the value of an Eisenstein series if $\rho$ is a complex multiplication, and otherwise a generalization of it.

\aster

When $\Omega$ is a lattice and $\rho \, \Omega \subset \Omega$, $\wp_{\rho}(z)$ coincides with the classical Weierstra{\ss}  function $\wp(z | \Omega)$. There are no convergence issues.   With the restrictions we have placed on $\rho$ (including  $\Delta =\{ 0, 1  \}$) there are, as previously asserted, just three cases. In each case, a fundamental domain for $\Omega \subset \C$ is the remainder set $R= \{ \sum_{k=1}^{\infty} \delta_k / \rho^k, \,  \delta_k \in \Delta \} $ and the collection $ \{ \omega + R : \omega \in \Omega \}$ is a periodic  tiling of $\C$ by copies of the tile $R$, and the corresponding $\wp_{\rho}(z)$ is an even function of $z$.

The new cases, of course, will be those for which $\rho$ is not a quadratic imaginary integer of norm 2; the $\rho$-integers do not form a lattice; and the associated inflationary tiling is not periodic. First let us look at the three quadratic imaginary  lattice examples in more detail.

\subsubsection{$\rho = i \sqrt{2}$}

For $\rho = i \sqrt{2}$ the remainder set is the rectangle
$$R = \left\{ x+i y : 0 \leq x \leq 1, \, 0 \leq y \leq \sqrt{2} \right\};$$ thus $R$ is the closure of a fundamental domain for the lattice $\Omega$ generated by $\{ 1, i \sqrt{2} \}$.

This positional representation and  $R$ can be expressed in terms of an integer radix because  $z = \sum \delta_k\, (i \sqrt{2})^{-k}$ can be rewritten as
\bean
 z &=&  \sum \left( \delta_{4k} (-2)^{-2k} + \delta_{4k+2} (-2)^{-(2k+1)} \right)   + \\
  	& & \qquad      i \sqrt{2} \sum \left( \delta_{4k+1} (-2)^{-2k} +  \delta_{4k+3} (-2)^{-(2k+1)} \right)  \\
  &=& x +i \sqrt{2} \, y ,
\eean
where $x, y \in \R$ are arbitrary numbers expressed in radix $-2$.
The units in $\Q(i\sqrt{2})$ are $\pm1$. The above argument shows that the $\rho$-integers are the same for radices $\pm i \sqrt{2}$.

 The $\rho$-duplication formula eq(\ref{rho wp dup}) specializes to the  complex multiplication formula
 \be
-2 \,  \wp(i \sqrt{2}\,  z  \, : \Omega) = \wp(z  : \Omega ) + \wp \left(z - i\sqrt{2}/2 \,  : \Omega \right) - \wp \left( - i\sqrt{2}/2   \, : \Omega \right).
 \ee
 
The numbers $\omega_1= -1$ and $\omega_2=-i\sqrt{2}$ are also lattice generators, and their sum is $\omega_3=-1-i \sqrt{2}$.  Choosing  $z=-1/2$, we obtain an identity for  the $\wp$-function evaluated at half periods,\footnote{We drop $\Omega$ to simplify the notation.} 
$$ 
-2 \wp ( \omega_2/2 )= \wp (\omega_1/2) + \wp(\omega_3/2)  - \wp(\omega_2/2) ,
$$
that is,
$$  \wp (\omega_1/2)  +\wp(\omega_2/2)  + \wp(\omega_3/2) = 0 . $$
This is a standard, and fundamental,  identity for the roots of $\wp'(z)$.

\subsubsection{$\rho = \frac{1 + i\sqrt{7}}{2}$}

$\rho = \frac{1 + i\sqrt{7}}{2}$ can be treated similarly. 

\subsubsection{$\rho = -1 +i $ \, : the lemniscatic case}

The lattice generated by 1 and $1+i$ is the ring $\Q[i]$ of Gaussian integers;  $1$ and $-1+i$ also generate $\Q[i]$, so the theory of elliptic functions is the same.  From the perspective of the associated positional representation, there is a difference. $\Delta[1+i]$ is not the full ring $\Q[i]$. One finds $\Delta[1+i] \cup(i - \Delta[1+i]) =\Q[i]$ and $\Delta[1+i] \cap (-\Delta[1+i]) = \emptyset $, whereas $ \Omega = \Delta[-1+i] =\Q[i]$.  In both cases, the remainder set is the so-called ``twindragon"  fractal region; it tiles $\Delta[1+i]$ and also $\C$. 

\begin{figure}[h]
\begin{center}
\caption{Decomposition of the Gaussian integers $\Q[i]$ (in a square neighborhood centered on the origin) as the disjoint union $\Delta[\rho] \cup (i - \Delta[\rho] )$ for radix $\rho = 1+i$. The blue points illustrate $\Delta[\rho]$.}
\label{rho=1+i_Omega}
\includegraphics[width=1.75in]{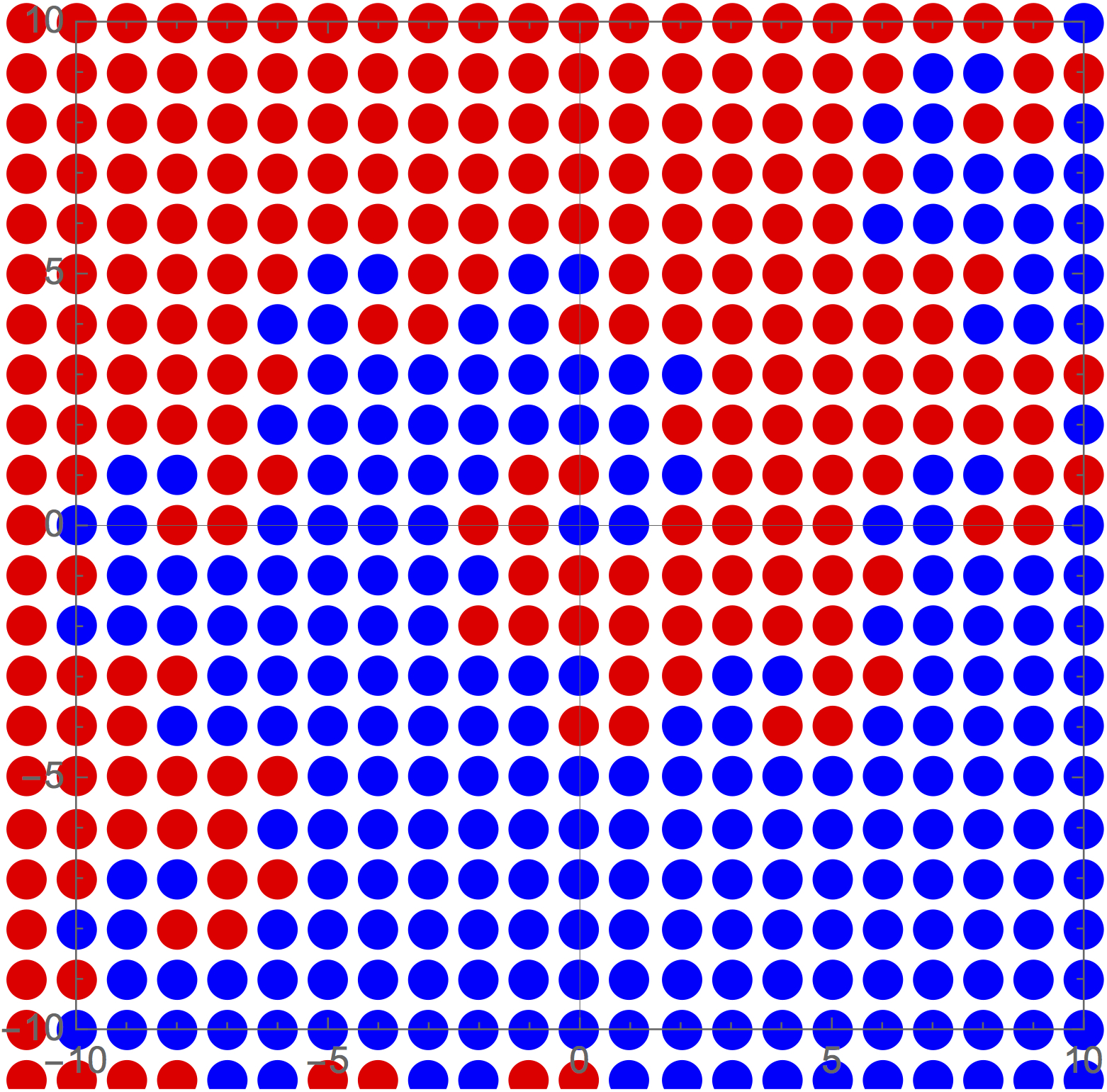}
\end{center}
\end{figure}

At this point it will be convenient to make a  change of notation. In the defining eq(\ref{rho-p-function})  for $\wp_{\rho}(z) $ we wrote $(z+\omega)^2$ where it is customary to write $(z-\omega)^2$ in order to preserve compatibility with the classical notation for the Hurwitz zeta function, and the $\rho$-generalization of it. This notation has the consequence that the poles of $\wp_{\rho}(z)$ do not in general coincide with the set of $\rho$-integers. It made no difference in the previous discussion of $\wp(z)$ because  $ - \Q[i] =\Q[i]$. But now, and for our later work,  it makes a difference because we will be considering aperiodic tessellations, so it will be preferable to define
\be
\wp_{\rho}(z) := \frac{1}{z^2} + \sum_{\Omega^*} \left( \frac{1}{ (z-\omega)^2} - \frac{1}{\omega^2} \right).
\label{wp_rho-def2}
\ee
In this notation, the $\rho$-duplication formula becomes
\be
\rho^2 \, \wp_{\rho}(\rho z) =   \wp_{\rho}( z) + \wp_{\rho}( z-1/\rho) - \wp_{\rho}( -1/\rho) .
\label{wp_rho-dup2}
\ee

For $\rho = 1+i$ the series for $\wp_{\rho}(z)$ is a 2-dimensional analogue of the series for the polygamma function $\psi^{(0)}(z)$, which runs over the non-negative integers. $\psi^{(0)}(z)$ is a kind of `one-sided' cotangent function (cp. eq(\ref{zeta_rho(1,z)})). Similarly, $\wp_{1+i}(z)$ is a `one-sided' Weierstra{\ss} $\wp$-function. The $(1+i)$-duplication formula is
\be
2i \, \wp_{1+i}((1+i)z) = \wp_{1+i}(z) +  \wp_{1+i} \left( z -\frac{1-i}{2} \right) -  \wp_{1+i} \left( -\frac{1-i}{2} \right) .
\label{1+i-dup}
\ee

\aster

With $\rho = (-1+i) $ we return to the realm of ordinary elliptic function theory and complex multiplication. The set of $\rho$-integers is the lattice  $\Delta[-1+i] = \Q[i]$, and  $\wp_{(-1+i)}(z) = \wp(z : \Q[i]) =\wp(z) $.\footnote{Again we drop the lattice $\Q[i]$ from the notation for simplicity.}

Observe that $1/\rho = \frac{1-i}{2}$  is a zero of $\wp^{\prime}_{1+i}(z)$ (it is congruent to  a half-period).  The $\rho$-duplication formula reduces to
\be
-2 i \, \wp ((-1+i)z ) = \wp(z ) + \wp \left( z+ \frac{1+i}{2}    \right) - \wp \left(  \frac{1+i}{2}   \right) .
\label{-1+i complex mult}
\ee
This is easily transformed into a formula with  complex multiplication $(1+i)$ for $\wp(z)$. Substitute $z = i t$ to find
$$
(-2i) \wp((1+i) t ) = 
	\wp( i t ) + \wp \left( i t + \frac{1+i}{2}   \right) - \wp \left( \frac{1+i}{2}   \right) ,
$$
and use the defining series to show that $\wp (i t ) = -\wp( t )$. Hence   eq(\ref{-1+i complex mult}) becomes
\be
2i \,  \wp((1+i) z ) = 
	  \wp(  z ) + \wp \left(  z + \frac{-1+i}{2}   \right)  + \wp \left( \frac{1+i}{2}   \right) .
\label{1+i complex mult raw}
\ee
Now substitute $z=1/2$ to find
$$
2i \, \wp \left( \frac{1+i}{2} \right) = \wp \left(   \frac{1}{2} \right) + \wp \left(  \frac{i}{2}   \right) +\wp \left( \frac{1+i}{2} \right),
$$
that is,
$$ 2 i \,  e_3=  e_1 + e_2 + e_3 = 0 ,$$
whence $\wp \left( \frac{1+i}{2} \right) = 0$ and eq(\ref{1+i complex mult raw}) assumes its final form,
\be
2i \,  \wp((1+i) z ) = 
	  \wp(  z ) + \wp \left(  z + \frac{-1+i}{2}   \right)   .
\label{1+i complex mult}
\ee
This brings us back to Fagnano's  complex multiplication formula, which can be obtained from eq(\ref{1+i complex mult}) by applying the addition formula for $\wp(z)$ to the second term on the right, and recalling that $\wp^{\prime}((-1+i)/2)=\wp((-1+i)/2)=0$.  We immediately find Fagnano's formula eq(\ref{Fagnano complex mult}):
$$
\wp((1+i)z) = -\frac{i}{8} \left( \frac{\wp^{\prime}(z) }{ \wp(z) } \right)^2 .
$$

\section{Examples with non-unique representation of $\rho$-integers}
\label{sec: non-unique examples}

Let $p$ and $q$ be polynomials in the indeterminate $X$ with coefficients in $\Delta$, and suppose that $p \neq q$. $p(\rho)$ and $q(\rho)$ are $\rho$-integers. If $p(\rho)=q(\rho)$ then the two polynomials represent the same $\rho$-integer, and this equality is an algebraic equation for $\rho$. This shows that non-uniqueness of the representation can only occur for certain algebraic values of the radix. Moreover, every pair of representations of the same $\rho$-integer leads to an algebraic equation which has $\rho$ as a root.

In our case, $1 < | \rho | \leq 2$ implies that the simplest examples arise as roots of the equations 
$$ x^2 =  \pm \, x + 1 ,$$
so that $\rho =  \pm \phi$.  Using positional notation relative to the radix,  the first case implies the identity $100 \cdot = 011 \cdot$; the second, $110  \cdot= 001 \cdot$, where `$\cdot$' denotes the radix point. Each equation is the underlying constituent of every non-unique representation of an integer for the associated radix. It can be used to construct a unique normal form for each integer. This construction will be taken up in the next sections. 

\label{phi-unique}
\subsection{Dimension 1: Properties of radix $ \phi = \frac{1+\sqrt{5}}{2} $}
\label{sec:radix phi}

Here is a summary of what we know about this radix so far: The  golden number $\phi =\frac{1+\sqrt{5}}{2}  \simeq 1.61$ is the larger root of $x^2=x+1$. The digits for the positional representation with radix $\phi$ are $\Delta = \{ 0, 1 \}$. The $\phi$-integers are the polynomials in $\phi$ with coefficients from $\Delta$, that is, $\Delta[\phi]$. Let $X$ denote an indeterminate. There are $2^n$ polynomials of degree less than $ n$ in $\Delta[X]$ but the specialization $X \rightarrow \phi$ reduces this number because of the identity $\phi^2 = \phi + 1$. In positional representation, this identity is $100 = 011$.\footnote{It is conventional to drop leading zeros, as well as the radix point when denoting integers.} A simple consequence is
 \be
 \phi^{2n} = 1+\sum_{k=1}^{n} \phi^{2k-1}, \quad n \in \Z^+ .
 \ee

There are two particular expressions that are useful `normal forms': in  the {\it maximal normal form}, the sequence of bits does not contain successive instances of  `1'; in the {\it minimal normal form}, the sequence of bits does not contain successive instances of `0'. For instance, the left side of $1001 =0111 = 111$ is in maximal normal form; the right side, in minimal normal form. The number of bits in the maximal normal form is  at least as large as in the minimal normal form.

The representation of  a $\phi$-integer  in either normal form is unique.   We shall use maximal normal form to construct a list of all  $\phi$-integers without repetition.

Every binary sequence that expresses an element in $\Delta[\phi]$  can be written in maximal normal form. Suppose $\delta = \delta_{n-1} \dots \delta_0$ is an arbitrary sequence of $n$ bits with  $\delta_{n-1} =1$. Suppose that $\delta_k \delta_{k-1}$ is the leftmost pair of successive `1's. If $k=n-1$ then  $\delta$ is equivalent to  $100\delta_{n-2} \dots \delta_0$. If $k<n-1$ then there is a 0 to the left of $\delta_k$ so the identity enables the replacement $0 \delta_k \delta_{k-1} \rightarrow 100$. If this introduces a sequence of successive `1's, we apply the first case. Repetition of the argument, i.e. induction,  completes the conversion to maximal normal form. 

Suppose that $\Omega_n$ denotes the set of all $n$-digit  binary sequences in maximal normal form for radix $\phi$. Set $\Omega_1 = \{ 0, 1 \}, \Omega_2 = \{ 0,1,\phi \}$. 

Recall that $\nabla_n$ denotes the set of $n$-digit $\phi$-integers whose leading digit is 1. For $\rho=\phi$ one finds the disjoint union
\be
\left.
\ba{rcl}
\nabla_1 &=& \{ 1 \} \\
\nabla_2 &=& \{ \phi \} \\
\nabla_n &=& (\phi\, \nabla_{n-1}) \cup (1 + \phi^2\, \nabla_{n-1})
\label{nabla_n for phi } 
\ea
\right\}
\ee
The number of elements in $\nabla_n$ is the Fibonacci number $F_n$.

From $\Omega^* = \bigcup_{n \geq 1} \nabla_n$ find
\bean
\Omega^* &=& \bigcup_{n \geq 1} \nabla_n \nonumber \\
	&=&  \bigcup_{n \geq 1} ( \phi \, \nabla_{n-1} 
		\cup  (1 + \phi^2\, \nabla_{n-1}) \nonumber \\
	&=&( \phi \, \Omega^*_{n-1} ) \cup (1 + \phi^2 \, \Omega^*_{n-2})
\eean 
so
\be
\Omega_n = ( \phi \, \Omega_{n-1} )  \cup (1 + \phi^2 \, \Omega_{n-2}) .
 \label{phi omega_n recursion}
\ee
The density of $\Omega$ in $\Z$ is $\lim_{n \rightarrow \infty} F_{n+2}/2^n =0$. 

\aster

Let $m$ denote Lebesgue measure on the plane induced by the Euclidean metric. Recall that two measurable subsets of $\C$ are {\it essentially disjoint} if the measure of their intersection is 0, and they are {\it essentially identical} if the measure of their intersection is equal to the measure of each of the sets. Thus  $m(A \cup B) = m(A) + m(B)$ if and only if $A$ and $B$ are essentially disjoint. 

In general, the remainder set $R$ for a positional representation with digits $\Delta=\{ 0,1 \}$  satisfies the recurrence
$$ \rho R = R \cup (1+R) .$$
If the sets on the right side are not essentially disjoint, then they will contain a set $S$ of positive measure. Now inflate $S$: $S \rightarrow \rho^n S$. The set of $n$-digit $\rho$-integers is the subset of $\ \rho^{n} R$ whose remainders are 0. For large enough $n$, $\rho^n S$ will contain a $\rho$-integer that, by its construction, will have more than one $\rho$-positional representation. This is the origin of the problem of non-unique representations of $\rho$-integers.

This description of the problem also contains the germ of a systematic alternative approach to finding unique representations for $\rho$-integers.

By introducing partial remainder sets that cover $R$ but are essentially disjoint,\footnote{These partial remainder sets are the `tiles' for the tessellation defined by the positional representation.} overlaps will have  measure zero. Let us follow this approach for $\rho=\phi$. The remainder set is $R = [0,\phi]$. Introduce $R_1= [0,1]$ and $R_2 =[0,\phi-1]$. Then
\bean
\phi \, R_1 &=& R_1 \cup (1 + R_2) \\
\phi \, R_2 &=& R_1 .
\eean
Letting $\Omega_{j,n}$ denote the set of $n$-digit $\phi$-integers in $\rho^n R_j$, these equations imply
\be
\left.
\ba{rcl}
\Omega_{1,n} &=& \Omega_{1,n-1} \cup \left( \phi^{n-1} + \Omega_{2,n-1}  \right) \\
\Omega_{2,n} &=&  \Omega_{1,n-1}
\ea
\right\}
\ee
whence
\be
\Omega_{1,n} =  \Omega_{1,n-1} \cup \left( \phi^{n-1} + \Omega_{1,n-2} \right), 
\label{phi recursion 2}
\ee
a relation similar to the last line of eq(\ref{phi omega_n recursion}). In fact, both recursions produce the same set of $\phi$-integers: eq(\ref{phi omega_n recursion}) constructs integers by inserting low-order digits; eq(\ref{phi recursion 2}) by inserting high order digits. That is, in any positional representation $x_{n-1} \dots x_0$ for a $\phi$-integer,  recursively make the substitution $ 0 1 1 \rightarrow 100$ so that the sequence $11$ no longer occurs.

\aster

The quadratic equation for $\phi$ implies that elements of $\Omega$ can be written as linear combinations of $1$ and $\phi$.  Explicitly, $m  \phi + n \in \Omega$ if and only if $m, n \in \Z^+ \cup \{0\}$ and $n$ satisfies the condition
$$  \left(  
	n =  
	 \lfloor m/\phi \rfloor  \right)  \vee  \left(   n = 1+ \lfloor m/\phi \rfloor \right)  		\vee 
	\left(  (n = 2+ \lfloor m/\phi \rfloor) \wedge   (2+ \lfloor m/\phi \rfloor > \phi + 	2+ \lfloor (m+1)/\phi \rfloor ) 
\right),   
$$
where $\lfloor x \rfloor$ denotes the least integer greater than or equal to $x$. This condition implies that for a given $m$, $n$ can always assume the first two of the  three values 
$\{   \lfloor m/\phi \rfloor   , 1+ \lfloor m/\phi \rfloor  , 2+ \lfloor m/\phi \rfloor   \}$ 
but not any value that is not one of the three. This provides a simple way to bound the possible values of $m \phi +n$. If $\omega = m \phi +n \in \Omega$, then
\be
2m/ \phi  <  m   \phi+ \lfloor m/\phi \rfloor   \leq \omega  \leq 
m  \phi +  m/\phi +2 = m (2 \phi-1)   +2.
\ee

\subsection{Dimension 1: Properties of radix $-\phi$}
\label{subsec:radix -phi}

Now let $\rho = -\phi$. Some integers have infinitely many representations and, more generally, $( -\phi)$-integers need not have a maximal representation.  Indeed, $\rho$ satisfies $x^2 +x =1$ which is equivalent to the positional identity $110 = 001$. It follows that $1$ has infinitely many representations, viz.
$$ 1\cdot  = 110 \cdot = 11010 \cdot =  \dots $$

Now consider the minimal representation, constructed by the substitution $110 \rightarrow 001$. The integer $\rho +1$ has the  representation $11 \cdot 0 = 00 \cdot 1$, which is a remainder. It is the only $\rho$-integer that is also a remainder.  We amend the definition to exclude polynomials $p(\rho)$ that are equal to a remainder. Suppose  $p, q, r \in \Delta[X]$ are polynomials and the degree of $q$ is $m$.  Furthermore suppose that the radix $\rho$ satisfies $ X^m p(X) = q(X)$. Then the polynomial $p(\rho) = q(\rho)/\rho^m$ is a remainder, hence not a $\rho$-integer. This is the situation we encounter with $\rho = -\phi$.

\aster 
The $(-\phi)$-integers are signed. Denote the set of $(-\phi)$-integers by $\Delta[\rho]$ and the subset of $n$-digit $(-\phi)$-integers by $\Delta_n[\rho]$.  We have
\be
\omega = \sum_{k \geq 0} \delta_k (-\phi)^k = \sum_{k \geq 0} \delta_{2k} \phi^{2k} - \phi \sum_{k \geq 0} \delta_{2k+1} \phi^{2k} 
\ee
where each sum is finite. Thus   $\Delta[\rho]$ is neither bounded below nor above. 

\aster
The construction of a complete set of unique representatives of the $\rho$-integers proceeds by using the equivalence of digit sequences $110 = 001$ to reduce representations to the least significant digit, i.e., to lower degree. This process is facilitated by introduction of {\it even} and {\it odd} $\rho$-integers. We shall say that a representation is {\it even}, resp. {\it odd}, if the least significant digit is 0, resp., 1. This notion depends on the form, but not the numerical value, of the representation. For instance, $\rho^2+ \rho = 110$ is even, but $1 = 001$ is odd.

Employing reduction to minimal normal form, set $\nabla^0_1 = \{ 0 \}$ and, for $n>1$, let $\nabla^0_n$ denote a unique set of even $n$-digit $\rho$-integers  with leading digit 1, and $\nabla^1_n$ the corresponding set of odd $\rho$-integers. $\nabla^1_2$ requires special attention because $11\cdot 0= 00 \cdot 1 $ is a remainder.  Hence $\rho + 1$ `looks like' a $\rho$-integer but is not one, so $\nabla^1_2 = \emptyset$. The following recursion is a constructive realization of a complete set of unique representatives for the $(-\phi)$-integers. For $n \leq 3$,
\be
\left.
\ba{rcl}
\nabla^0_1 &=& \{ 0 \}, \qquad  \nabla^1_1 = \{ 1 \} \\
\nabla^0_2  &=& \{ \rho \} , \qquad  \nabla^1_2 = \emptyset \\
\nabla^0_3  &=& \{ \rho^2 \} , \,\,\, \,  \quad \nabla^1_3 = \{ \rho^2 +1,  \rho^2 +\rho  + 1 \} 
\ea
\right\},
\label{-phi V_n initial cond}
\ee
whereas for $n>3$, 
\be
\left.
\ba{rcl}
\nabla_n &=& \nabla^0_n \cup \nabla^1_n \\
\nabla^0_n &=& (\rho \, \nabla^0_{n-1}) \cup  (\rho+\rho^2 \, \nabla^0_{n-2})  \\
\nabla^1_n &=&    (1 +\rho \, \nabla^0_{n-1})  \cup  (1 +\rho \, \nabla^0_{n-1}) 
\ea
\right\}
\label{-phi V_n recursion}
\ee
As usual, $\Omega^*_n = \bigcup_{k=1}^n \nabla_n $ and $\Omega = \{ 0 \} \cup \Omega^*$. Set $\Omega^0 =  \Omega \cap \left(  \bigcup_{k=2}^n \nabla^0_n \right) $ -- these are the even $\rho$-integers with at most $n$-digits.

The construction of  $\nabla^1_n$ is intertwined with the $\nabla^0_k$.  According to the second to last line of eq(\ref{-phi V_n recursion}),  the even representations can be constructed without reference to the odd ones. One consequence is that this set of even $\rho$-integers can be used to define functions that have $\rho$-duplication formulae.  It follows that, with the initial conditions
$$
\Omega^0_1 =\{ 0 \}, \quad \Omega^0_1 =\{ 0, \rho \}, \quad \Omega^0_3 =\{ 0,\rho, \rho^2   \},
$$
the recursion is
\be
\Omega^0_n = ( \rho \,  \Omega^0_{n-1}) \cup (\rho + \rho^2 \, \Omega^0_{n-2} ).
\label{-phi even recursion}
\ee

The recursion summands are pairwise disjoint. Hence the number of elements in these sets can be calculated  by observing, for example, that $\#(\nabla^0_n) = \#(\nabla^0_{n-1}) + \#(\nabla^0_{n-2})$, and solving the recurrence with the aid of the initial values to find $\#(\nabla^0_n) = F_{n-1}$, and then  $\#(\Omega^0_n) =   F_{n+1}$.

\subsection{Dimension 2: Properties of radix $ i \sqrt{\phi} $}
\label{subsec: rho = i sqrt phi}

Now consider the radix $\rho = i \sqrt{\phi}$ with remainders $\Delta=\{0,1 \}$ as above.  We shall obtain a unique representative for each $\rho$-integer by expressing the $\rho$-integers  in terms of the $\phi$-integers, and then applying the unique expression for those integers realized by the set $\Omega$ to obtain a complete set of unique representatives for $\Delta[ \rho ]$.

The remainder set $R =\{ x : x = \sum_{k=1}^{\infty} \delta_k/\rho^k  \} ,  \delta_k \in \Delta, $ is a rectangle in $\C$ whose vertices are
 $$ \sum_{k=1}^{\infty} \frac{1}{\rho^{4k}}, \quad  \sum_{k=1}^{\infty} \frac{1}{\rho^{4k+1}}, \quad \sum_{k=1}^{\infty} \frac{1}{\rho^{4k+2}}, \quad \sum_{k=1}^{\infty} \frac{1}{\rho^{4k+3}}, $$ 
 that is, $\phi - 1, -i \sqrt{\phi} (2-\phi) , (2-\phi), i \sqrt{\phi}(2-\phi)$.
 
An arbitrary polynomial $\omega = \sum_{k=0}^{n-1} \delta_k \rho^k $ of degree $n-1$   in $\Delta[ \rho ] $ can be expressed as
$$ 
  \sum_{k=0}^{n-1} \delta_{4k} \rho^{4k} +  
 	\sum_{k=0}^{n-1} \delta_{4k+1} \rho^{4k+1}  +
		  \sum_{k=0}^{n-1}  \delta_{4k+2} \rho^{4k+2}  +
		  	 \sum_{k=0}^{n-1}  \delta_{4k+3} \rho^{4k+3} . 
$$
 Then
\bea
\omega &=& 
 	\sum_{k=1}  \delta_{4k} \phi^{2k}  +
		i \sqrt{\phi} \sum_{k=1}  \delta_{4k+1} \phi^{2k}   +
			- \phi \sum_{k=1}  \delta_{4k+2} \phi^{2k}  +
			i \phi  \sqrt{\phi}	\sum_{k=1}  \delta_{4k+3} \phi^{2k} \nonumber  \\
	&=& \omega_0 +
			(i \sqrt{\phi}) \, \omega_1 -
				 \phi \,  \omega_2   -
				  	(i \sqrt{\phi}) \, \phi \,\omega_3 \nonumber  \\
	&=& \left( \omega_0 -\phi \,\omega_2  \right) + (i \sqrt{\phi})\left( \omega_1 -\phi \, \omega_3 \right).
\label{eq:i rt phi integers}
\eea
If the $\omega_k - \phi \, \omega_{k+2}$ are restricted to a complete set of unique $(-\phi)$-integers,  then eq(\ref{eq:i rt phi integers}) yields a unique expression for each $\rho$-integer in terms of powers of the radix $\rho$ and $(-\phi)$-integers. That this expression for the $\rho$-integers still requires only 2 digits for a complete set of positional representations relative to radix $\rho^2$  is a consequence of $| \rho^2 | =| -\phi | \leq 2$.
 	
\section{Identities for partial sums of meromorphic functions}
\label{numerical exact identities}

In this section we consider partial sums of the $\rho$-Hurwitz $\zeta$-function and show that they satisfy certain identities even when the complete series is divergent. We begin with the general case, for which representations of $\rho$-integers are unique. This  is the simplest, but also the least interesting, because it is not related to tessellations.

\subsection{If $\rho$-integer representations are unique  } 

Let $\Omega_n = \Delta_n[\rho]$ denote the set of $n$-digit $\rho$-integers, and assume that the representation is unique. Thus the  coefficient of $\rho^{n-1}$ is 1.  If representations of  $\rho$-integers are unique, the $2^n$ $n$-digit polynomials in $\rho$ represent distinct integers and one has the identity
$$ \Omega_n =  (\rho \, \Omega_{n-1} ) \cup ( 1 + \rho \, \Omega_{n-1}) .$$
Introduce
\be
f_n(s,z) = \sum_{\Omega_n} \frac{1}{(z + \omega_n)^s}.
\ee
A simple calculation shows that
\be
\rho^s \, f_n (s,\rho z) = f_{n-1}(s,z) + f_{n-1}(s, z + 1/\rho).
\ee
This $\rho$-duplication formula for $f$ is an identity for finite sums so convergence is not an issue.

The set of all $\rho$-integers is
$$ \Omega = \lim_{n \rightarrow \infty} \Omega_n, $$
and the $\rho$-Hurwitz $\zeta$-function $\zeta_{\rho}(s,z)$ is defined as
$$ \zeta_{\rho}(s,z) = \sum_{\Omega} \frac{1}{(z+\omega)^s},$$
if the sum converges appropriately, whence, under the same conditions,
\be
\zeta_{\rho}(s,z) =\lim_{n \rightarrow \infty} f_n(s,z) .
\ee
Therefore  $\zeta_{\rho}(s,z)$ satisfies the same $\rho$-duplication formula, viz.,
\be
\rho^s \, \zeta_{\rho} (s,\rho z) =  \zeta_{\rho}(s,z) + \zeta_{\rho}(s, z + 1/\rho).
\ee

\subsection{If $\rho = \phi$}

Let $\rho = \phi$. Representations of $\rho$-integers are not unique, so we must re-examine identities for partial sums.  Introduce  
\bea
f_{n}(s,z) &=&  \sum_{\omega \in \Omega_n} \frac{1}{(z+ \omega)^s} , \\
f_{n}(z) &=& \frac{1}{z}+ \sum_{\omega \in \Omega^*_n} \left( \frac{1}{z+ \omega} - \frac{1}{\omega} \right) .\\
\eea
Since the sums are finite there is no question of convergence for any $s$. 
A calculation shows
\bea
\rho^s f_{n}(s,\rho z)&= & f_{n-1}(s,z) + \frac{1}{\rho^{s}} f_{n-2}(s,z/\rho + 		1/\rho^2),\\
\label{f_n(s,z)}
\rho f_{n}(\rho z)&= & f_{n-1}(z) + \frac{1}{\rho} f_{n-2}(z/\rho + 1/\rho^2) - 	\frac{1}{\rho} f_{n-2}( 1/\rho^2).
\label{f_n(z)}
\eea
For each $n \in \Z^+$, these are identities for finite sums.
 
Now consider the limiting cases.  Let $\R_{\Omega}$ denote the complement of the set of poles $-\Omega$ in $\R$. In section \ref{sec: R+ conv} we  show that 
\be
f_n(s,z) = \lim_{n \rightarrow \infty} f_{n}(s,z) = \sum_{\omega \in \Omega} \frac{1}{(z+\omega)^s}
\ee
exists and the series converges uniformly  on compact subsets of $\R_{\Omega}$ for $\Re(s)>1$, and that 
\be
f(z) = \lim_{n \rightarrow \infty} f_{n}(z) = \frac{1}{z} +  \sum_{\omega \in \Omega_n} \left( \frac{1}{z+ \omega} - \frac{1}{\omega} \right)
\ee
also converges uniformly on compact subsets of $\R_{\Omega}$. As a consequence one sees that $f(z) = \cot_{\rho}(z)$ so
 
\bea
\rho^s f(s,\rho z)&= & f(s,z) + \frac{1}{\rho^{s}} f(s,z/\rho + 1/\rho^2) ,    \\
\rho  \cot_{\rho}(\rho z)&= & f(z) + \frac{1}{\rho}  \cot_{\rho}(z/\rho + 1/\rho^2) - \frac{1}{\rho}  \cot_{\rho}( 1/\rho^2) .
\label{n digit identity}
\eea
The $f(s,z)$ are proportional to the derivatives of $ \cot_{\rho}(z)$.

\subsection{If $\rho = -\phi$}
\label{subsec: rho=-phi}

We shall give an exact $\rho$-duplication formula for the subset $\Omega^0$ of even $(-\phi)$-integers.  The argument follows the same lines as the previous section but is based on  the defining recursion eq(\ref{-phi V_n recursion}) for the $\Omega^0_n$.  Put $(\Omega^0_n)^* = \Omega^0_n - \{ 0 \}$, $\Omega^0 = \lim_{n \rightarrow \infty} \Omega^0_n$, and $(\Omega^0)^* = \Omega^0 - \{ 0 \}$.

Introduce the $\rho$-Hurwitz $\zeta$-function by
\be
\zeta^0_{\rho}(s,z) := \sum_{\Omega^0} \frac{1}{(z+ \omega)^s},
\label{even -phi-Hurwitz}
\ee
and the functions
\bean
f^0_n(s,z) &:=&   \sum_{\Omega^0_n} \frac{1}{(z+ \omega_n)^s} \\
f^0_n(z) &:=&   \frac{1}{z} +\sum_{(\Omega^0_n)^*} \left(  \frac{1}{(z+ \omega_n)} - \frac{1}{\omega_n}  \right).
\eean
Then 
\bean
\zeta^0_{\rho}(s,z) &=& \sum_{n=1}^{\infty}  f^0_n(s,z) \\
\zeta^0_{\rho}(z) &=& \sum_{n=1}^{\infty}  f^0_n(z) .
\eean
 The difference between the even $(-\phi)$-integers and the  $\phi$-integers appears from eq(\ref{-phi even recursion}), 
$$ 
\omega_n \in \{ \rho \, \omega_{n-1}, \rho + \rho^2 \,\omega_{n-2}  \}, 
$$
which implies
\be
\rho^s f^0_n(s,\rho z) = f^0_{n-1}(s,z) + \frac{1}{\rho^s} f^0_{s,n-2}(z/\rho + 1/\rho) 
\label{-phi even numerics}
\ee
is an exact equality, whence, if the series converges,
\be
\rho^s\, \zeta^0_{\rho}(s,\rho z) =  \zeta^0_{\rho}(s, z) + 
	  \frac{1}{\rho^s}\zeta^0_{\rho}(s,z/\rho + 1/\rho) 
\label{(-phi)-zeta even dup}
\ee
is a $\rho$-duplication formula.

This duplication formula contains an analogue of the formula $\sum_{k=1}^{\infty} \frac{1}{(2k+1)^s} = (2^s-1) \zeta(s)$. Taking the limit as $z \rightarrow 0$, one finds
\be
(\rho^s -1) \sum_{\left(\Omega^0 \right)^*} \frac{1}{\omega^s} =
	 \sum_{\Omega^0 } \frac{1}{(1 + \rho \, \omega)^s} .
\ee

\section{Analytic continuation of $\zeta_{\phi}(s)$}
\label{sec: zeta analytic continuation}

We have  shown how an identity for the Riemann $\zeta$-function derived from the positional representation for radix $\rho = 2$ provides an analytic continuation for $\zeta(s)$ to the entire complex plane  (eq(\ref{zeta_rho-Riemann Laurent series})).  Now we will follow that line of reasoning for radix $\rho = \phi$. Once again we emphasize that  positional representations for $\rho$-integers are not unique for this radix.  This is the interesting new fact.

Summarizing previous notation, let $\Omega_n$ denote the set of $\rho$-integers defined by the recursion
$$ 
\left.
\ba{rcl}
\Omega_1 &=& \{ 0,1 \} \\
\Omega_2 &=& \{ 0,1, \phi   \} \\
\Omega_n &=& (\phi \, \Omega_{n-1} ) \cup (1 + \phi^2 \, \Omega_{n-2} )
\ea
\right\} ,
$$
and put $\Omega = \lim_{n \rightarrow \infty} \Omega_n$. The set $\Omega_n$, resp. $\Omega$, is a complete set of unique representations for the $n$-digit, resp., all, $\rho$-integers. The number of elements in $\Omega_n$ is $F_{n+2}$.  Let $\Omega_n^*$, resp. $\Omega^*$, denote the subset of non zero elements of $\Omega_n$, resp. $\Omega$.

The $\phi$-Hurwitz $\zeta$-function is
$$ \zeta_{\phi}(s,z) = \sum_{\Omega} \frac{1}{ (z +  \omega)^s} .$$
Temporarily abbreviate $\zeta_n(s,z) = \sum_{\Omega_n} (z+ \omega_n)^{-s}$,  and calculate
\bean
 \zeta_{n}(s,z) 
 	&=& \sum_{\Omega_{n-1}} \frac{1}{ (z + \omega_{n-1})^s} +
 	\sum_{\Omega_{n-2}} \frac{1}{ (z + 1 + \phi^2 \,  \omega_{n-2})^s}  \\
	&=& \frac{1}{\phi^s}  \sum_{\Omega_{n-1}} \frac{1}{ (z/\phi + 		\omega_{n-1})^s} +
	\frac{1}{\phi^{2s}}  \sum_{\Omega_{n-2}} \frac{1}{ (z/\phi^2 + 		1/\phi^2 + 	 \omega_{n-2})^s}  \\
	&=& \frac{1}{\phi^{s}} \zeta_{n-1}(z/\phi + 1/\phi) + \frac{1}{\phi^{2s}} \zeta_{n-2}(s,z/\phi^2 + 1/\phi^2).
\eean
Thus  
\be
\phi^s \,  \zeta_{n}(s,\phi z)  = \zeta_{n-1}(s, z) + 
	 \frac{1}{\phi^s} \zeta_{n-2}(s, z/\phi + 1/\phi^2) ,
\label{phi zeta_n dup}
\ee
and the limit as $n \rightarrow \infty$ yields the $\phi$-duplication formula for $\zeta_{\phi}(s,z)$:
\be
\phi^s \,  \zeta_{\phi}(s,\phi z)  = \zeta_{\phi}(s, z) +  \frac{1}{\phi^s} \zeta_{\phi}(s, z/\phi + 1/\phi^2) .
\ee
From the defining equation for $\zeta_{\phi}(s,z)$, observe that $h_n(s,z):=\phi^s\, \zeta_n(s, \phi z) -  \zeta_{n-1}(s,  z) $ is holomorphic at $z=0$, so
\bean
\lim_{n \rightarrow \infty} \lim_{z \rightarrow \infty} h_n(s,z) 
	&=&	( \phi^s -1) \sum_{\Omega^*} \frac{1}{\omega^s} \\
	&=&  ( \phi^s -1) \, \zeta_{\phi}(s) ,
\eean
whence, from eq(\ref{phi zeta_n dup}), 
\be
\phi^s  ( \phi^s -1) \, \zeta_{\phi}(s) =  \zeta_{\phi}(s,1/\phi^2) .
\label{phi Hurwitz to Riemann}
\ee

Now we are in a position to calculate the Laurent expansion of $\zeta_{\phi}(s,z)$ at $z=0$.  The auxiliary function
$$ h(s,z) :=\zeta_{\phi}(s,z) -\frac{1}{z^s} = \sum_{\Omega^*} \frac{1}{(z + \omega)^s} $$
is holomorphic at $z=0$. Calculate
$$ h(s,0) = \zeta_{\phi}(s) $$
and
$$ d^k h(s,z)/dz^k = (-1)^k \frac{ \Gamma(s+k)}{\Gamma(s) \Gamma(k+1)} \zeta_{\phi}(s +k).
$$
The Laurent expansion is
$$
\zeta_{\phi}(s,z) = \frac{1}{z^s} + \zeta_{\phi}(s) + 
	\sum_{k=1}^{\infty}  (-1)^k \frac{ \Gamma(s+k)}{\Gamma(s) \Gamma(k+1) } \zeta_{\phi}(s +k) \, z^k .
$$
Evaluate at $z=1/\phi^2$ and use eq(\ref{phi Hurwitz to Riemann}) to find
\be
\zeta_{\phi}(z) = \frac{\phi^{2s}}{\phi^{2s}-\phi^s-1}
	\left( 1 + 	
	\sum_{k=1}^{\infty}  (-1)^k \frac{ \Gamma(s+k)}{\Gamma(s) 		\Gamma(k+1) } \frac{  \zeta_{\phi}(s +k)}{\phi^{2(s+k)}  }
	\right).
\label{phi zeta Laurent series}
\ee
This is the analytic continuation of $\zeta_{\phi}(z) $. We shall show below that the series converges absolutely and uniformly on compact subsets of the complement of the set of poles $-\Omega$ for $\Re(s)>1$.

Comparison of eq(\ref{phi zeta Laurent series}) with the corresponding Laurent expansion for the classical Riemann $\zeta$-function, eq(\ref{zeta_rho-Riemann Laurent series}) shows that they have the same general form. The factor $ \frac{\phi^{2s}}{\phi^{2s}-\phi^s-1}$  preceding the series has a pole at the root of the denominator, that is, where $x=\phi^s$ satisfies $x^2 = x+1$, i.e., at $\phi^s=\phi$ or $\phi^s=1/\phi$. Hence $s= 1 $ is a pole.  The pole of the corresponding factor in eq(\ref{zeta_rho-Riemann Laurent series}) is the root of $\rho^s =2$ with $\rho=2$, i.e., of $2^s=2$. In both cases, the algebraic equation for the radix determines the pole.

The convergence in a right half plane of the series  eq(\ref{phi zeta Laurent series}) is enough to show that it provides an analytic continuation for $\zeta_{\phi}(s)$ to the whole plane. The series is not, however, useful for calculating values of the function because it converges slowly. 

\section{A 1-dimensional aperiodic tessellation and the function  $\cot_{-\phi} (z)$}
\label{sec: cot0_-phi}

Now we shall turn to the connection of these functions to tessellations, and, in particular, to aperiodic tilings. First  we focus our attention on radix $\rho = -\phi$ and the corresponding 1-dimensional tiling. The $\rho$-integers are unbounded from below and above.  Recalling the discussion preceding  eq(\ref{-phi even recursion}), a set $\Omega$ of unique representatives of the even $\rho$-integers was defined by that equation. $\Omega_n$  satisfies the recursion
$$
 \Omega^0_n =  \left(\rho\, \Omega^0_{n-1} \right) \cup \left(\rho+ \rho^2 \, \Omega^0_{n-2} \right) .
$$

Now define a cotangent function for this situation by 
\be
\cot^0_{\rho}(z):= \frac{1}{z} + 
	\sum_{\left(\Omega_n \right)^*} 
	\left( 
	\frac{1}{z - \omega}  + \frac{1}{\omega}  
	\right) .
\ee
This choice of sign for $\omega$ insures that the poles of $\cot_{\rho} (z)$ coincide with the $\rho$-integers.

A duplication formula for $\cot_{\rho}(z)$ can be found by repeating the argument that led to eq(\ref{(-phi)-zeta even dup}), or by integration of that equation with $s=1$. The result is
\be
\rho \, \cot_{\rho}( \rho z) = \cot_{\rho}( z) + \frac{1}{\rho} \,\cot_{\rho}( z/\rho - 1/\rho^2) -  \frac{1}{\rho} \, \cot_{\rho}( -1/\rho^2) .
\ee

The  $\rho$-integers are the poles, and the differences between successive poles is either $\phi$ or $\phi + 1$. The function $\cot^0_{\rho} z$ is approximately periodic\footnote{The phrase ``approximately periodic" is used in an informal sense.} on nearby inter-polar segments of the same length. That is, when the graph of $\cot^0_{\rho} z$ on a vertical strip between two successive poles is shifted onto a congruent such strip, the graphs of the function are approximately equal. Moreover, the graphs on the strip of width $\phi$ are approximately equal to the graphs on the strips of width $\phi+1$ after they are simultaneously scaled by a fixed factor. The precise relationship is given by the $\rho$-duplication equation.

\begin{figure}[t]
\begin{center}
\caption{Graph of $\cot^0_{\rho}(x)$ for the subset of even integers for  radix $\rho = - \phi $, showing  the location (indicated by the red bars) of the poles at the  even $\rho$-integers. The successive intervals constitute an aperiodic tessellation of $\R$.}
\label{-phi cot}
\includegraphics[width=2.5in]{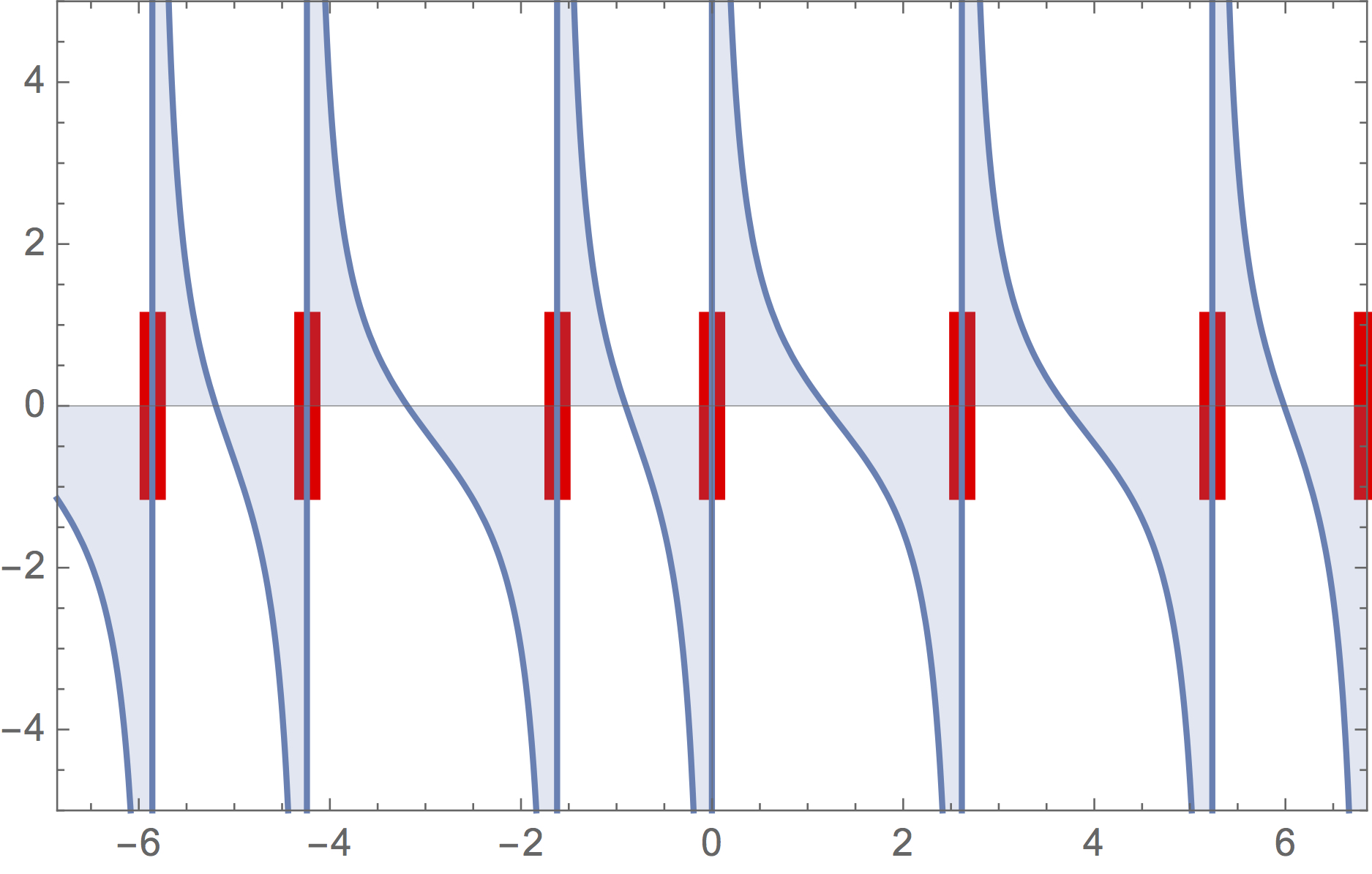}
\end{center}
\end{figure}

If the poles of $\cot^0_{\rho}(z)$ are interpreted as the endpoints of the intervals they delimit, these intervals are the tiles of an aperiodic  inflationary 1-dimensional tessellation. Hence there are two types of tiles: intervals of length $\phi$, and intervals of length $\phi +1$.  The poles, and therefore these tiles, are distributed aperiodically. Figure \ref{-phi cot} displays the graph of $\cot^0_{\rho}(x)$ and part of the aperiodic 1-dimensional tessellation formed by the even  $\rho$-integers.

\clearpage


\section{Function theory on inflationary tessellations.   The special tessellation}
\label{sec: function theory}

\subsection{Preparations}
\label{sec: preparations}

Now we  turn our attention to function theory on inflationary tessellations of the plane. 

Consider an aperiodic inflationary tessellation of $\C$ constructed from $n$ types of tiles defined with respect to the equivalence relation of congruence.  Since copies of a tile of a given type are congruent,  there is a 1-1 correspondence of their points.\footnote{Up to a set of measure 0.} 


The  special tessellation offers an interesting example. It is `special' because it is aperiodic; it consists of four tile types, all of which are rectangles; three of the four tile types are similar.   The radix (the `multiplier') is $\rho = i \sqrt{\phi}$ where $\phi = \frac{1+\sqrt{5}}{2}$ is the golden number; the tiles in the tessellation are rotated by the multiple of $ \pi/2$ that is the power of $\rho$  is associated with in the positional representation, but they are never reflected. This tiling  provides a convenient test laboratory for studying functions defined on an aperiodic tessellation.  

Consistent with previous notation,  $\Delta = \{ 0, 1 \}$ is the set of digits, and $ \Delta[\rho]$ denotes the set of polynomials in $\rho$ with coefficients in $\Delta$; these are the `integers' relative to radix $\rho$. There are repetitions. We shall construct a set of unique representations  $\Omega$ for a subset of the $\rho^2$-integers. Like the set of `even' $(-\phi)$-integers considered in section \ref{sec: cot0_-phi}, $\Omega$ is a set of unique representations, but it is not complete. Then  we will use $\Omega$ to construct a set of unique representations for a subset of the $\rho$-integers.

$\rho$ is the largest imaginary root of $x^4 + x^2 =1$. As we know (because $\rho^2=-\phi$), this implies that infinitely many elements of $\Delta[\rho]$ have more than one expression in $\Delta[\rho]$. For instance, 
\be
1 =\rho^{4m} + \sum_{k=1}^{m} \rho^{4k-2}, \qquad m \geq 1
\ee
Division by $\rho^2$ shows that a sum of $\rho$-integers can be a remainder: $1/\rho^2 = 1 + \rho^2$; indeed, $\rho^{-2} = -1/\phi =1- \phi  \sim -0.61$ has infinitely many representations as a sum of non-negative powers of $\rho$. The various realizations of $1+\rho^2 $ have the least absolute value of any element of $\Delta[\rho]$.

\aster
Now turn to the associated tessellation.   The remainder set is
\be
R =  \left\{ z: z = \sum_{k=1}^{\infty} z_k\, \rho^{-k} \right\}, \quad z_k \in \Delta  \, ,
\ee
that is, the rectangle
$$ R = \left\{ x+ i y : \rho^2/\phi   \leq  x \leq  1/\phi  , \,\,
				\rho^3/ \phi   \leq  y \leq \rho/\phi   \right\} . $$

An arbitrary point of $\C$ will be covered by a sequence of   inflations of  $R$  because  $0$ is an interior point of  $R$ so the $n$-th inflation contains the disk of radius $r^n$ centered $0$ if $r < 1/\phi$ because the disk of radius $r$ centered at $0$ is contained in $R$.  It follows that an essentially disjoint decomposition of $R$ that behaves well with respect to inflation by $\rho$ will lead to a tessellation of the plane. We will construct such a decomposition by the method introduced in the previous section. 

\begin{figure}[h]
\begin{center}
\caption{Radix $\rho = i \sqrt{\phi} $. Decomposition of the remainder set $R$ into an essentially disjoint union of remainder sets $R_1, R_2, R_3, R_4$. $R_1$ is the (magenta colored) rectangle in the first quadrant.}
\label{HLR_2D-rtPhi_R_decomp}
\includegraphics[width=0.8in]{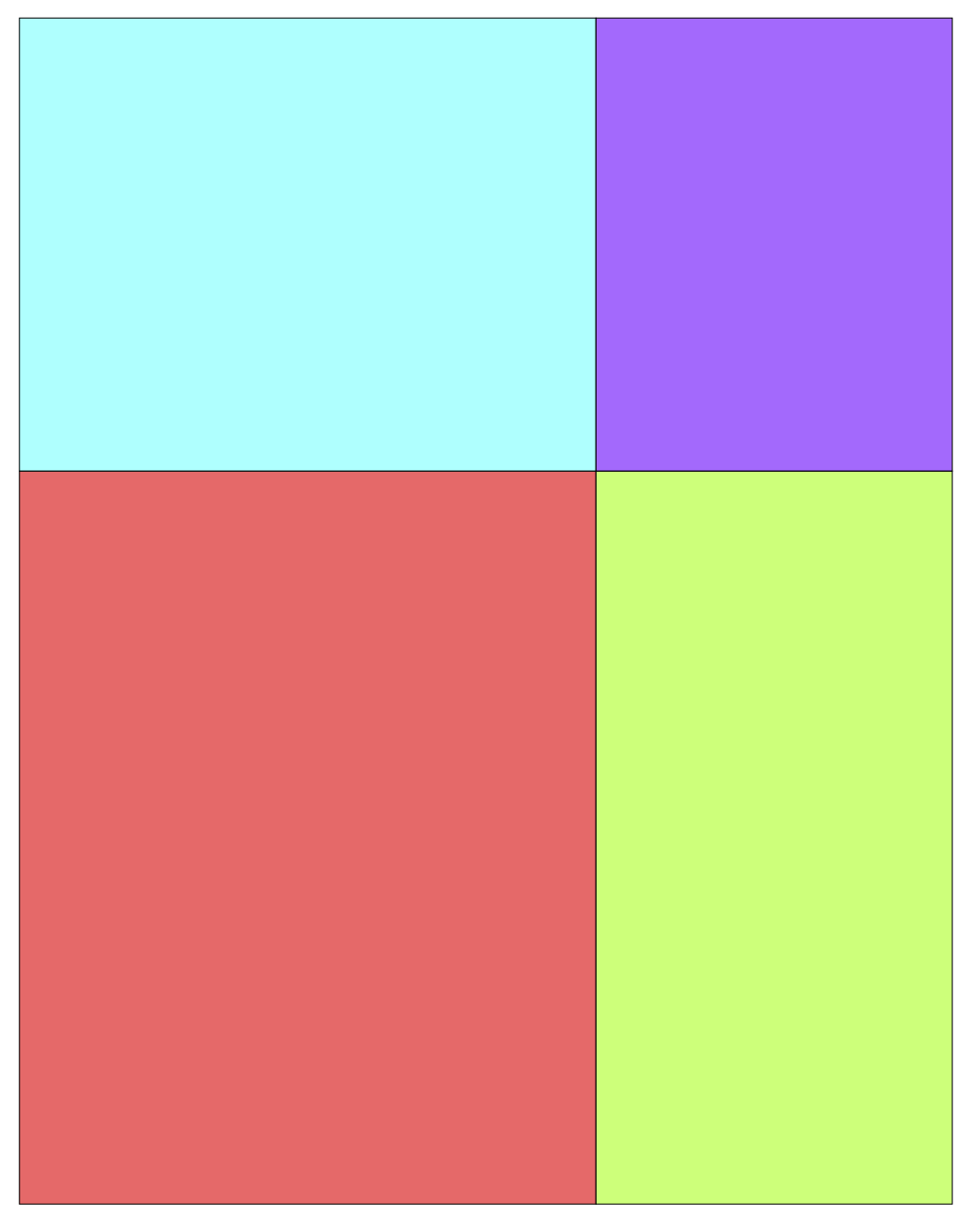}
\end{center}
\end{figure}

Although $R$ contains the origin as an interior point, its successive inflations cover  the plane but do not tile it because the decomposition  $R = (1+R) \cup R$  is not essentially disjoint. The translated copies overlap each other. Following the procedure previously used in 1-dimension \cite{HLR arXiv 3}.  $R$ can be written as the essentially disjoint union of 4 subsets $R_i$ that satisfy eq(\ref{rt-phi decomp}); cp. fig.\,\ref{HLR_2D-rtPhi_R_decomp}. The subsets are the parts of $R$ contained in the quadrants of the plane. This situation is described by the tessellation defined by  
\be
\left.
\ba{rcl}
\rho \,R_1 &=& R_2   \\
\rho \,R_2 &=&  R_3  \\
\rho \,R_3 &=& \left(1 + R_3 \right) \cup   R_4  \\
\rho \,R_4 &=& \left(1 + R_2 \right) \cup  R_1 
\ea
\right\}
\label{rt-phi decomp}
\ee

The tile types $R_1$, $R_2$ and $R_3$ are similar.  $R_1$ is the intersection of the remainder set $R$ with the first quadrant.  This  will be called the {\it special} tessellation.

\begin{figure}[h]
\begin{center}
\caption{Radix $\rho = i \sqrt{\phi} $. 10-digit tessellation of remainder set $R_1$. See text for details.}
\label{dim2 silver(0101) R1 10 digit}
\includegraphics[width=1.75in]{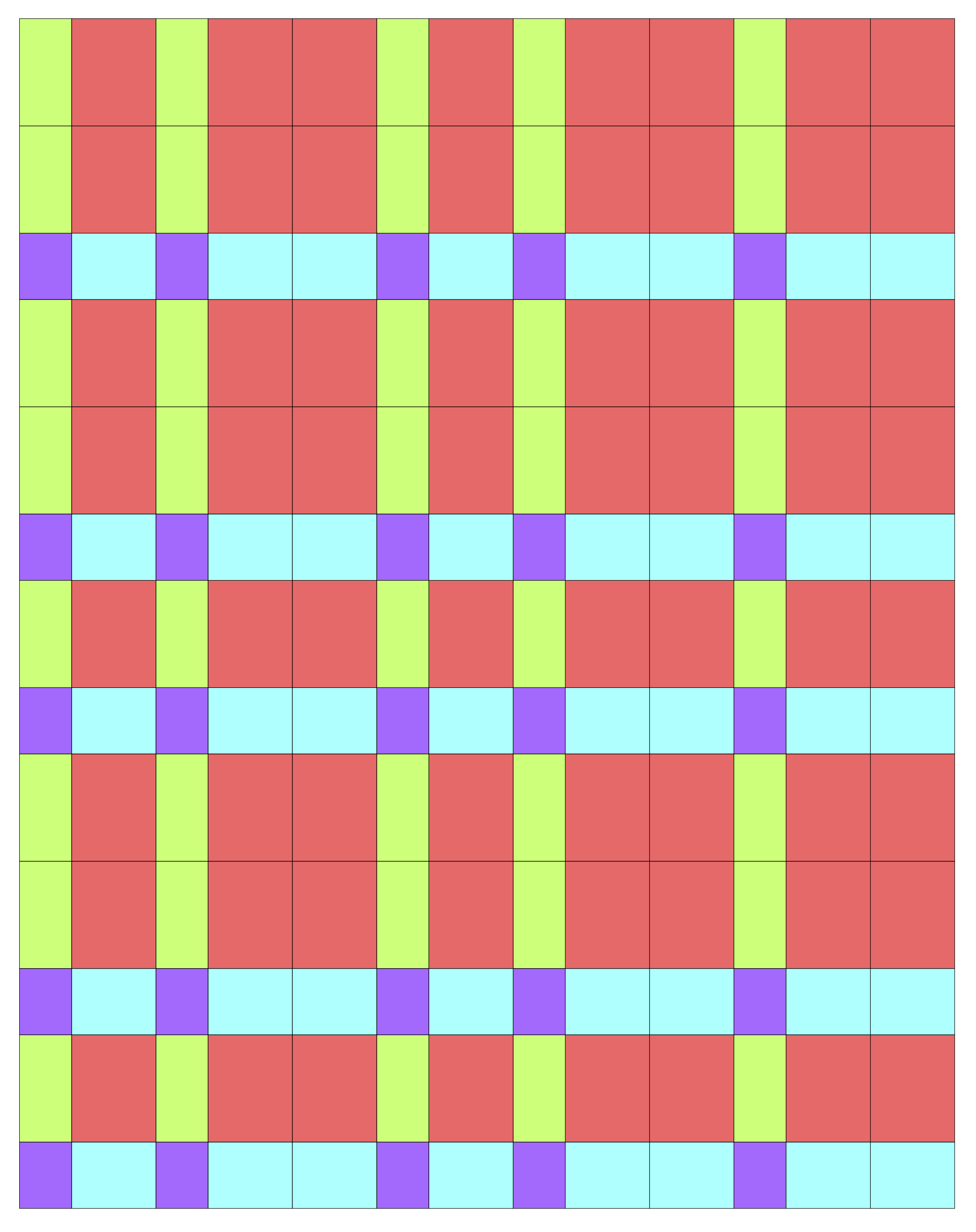}
\end{center}
\end{figure}

A portion of the 10-digit decomposition of $R_1$ is shown in fig.\,\ref{dim2 silver(0101) R1 10 digit}. Note that complex conjugation does not appear in the equations which implies that the tiles in the decomposition equations are not  reflected.  
\vspace{12pt}
\begin{it}
The special tessellation is aperiodic.\\
\end{it}

\, \vspace{12pt}
\pf The proof by irrationality\footnote{Cp. \cite{HLR arXiv 3}.} applies. The partition matrix, 
$$
U = \left(
\ba{cccc}
0 & 1 & 0 & 0 \\
0 & 0 & 1 & 0 \\
0 & 0 & 1 & 1\\
1 & 1 & 0 & 0
\ea
\right) \, ,
$$
has largest real eigenvalue $\phi = \frac{1+ \sqrt{5}}{2}$, so  $\rho = \sqrt{-\phi}$ is a multiplier (Note that $U$ is the partition matrix for a 2-dimensional tiling, so the eigenvalue is the square of the absolute value  of the  inflation factor. The tile type densities, which can be calculated from limits of quotients of the entries of $U^n$, are irrational, so the tiling cannot be periodic; cp. \cite{HLR arXiv 3}. $\Box$

\subsection{$\wp_{\rho}(z)$  and its $\rho$-duplication formula}
\label{subsec: wp_(-phi)}

According to eq(\ref{eq:i rt phi integers}), each  $\rho$-integer can be written in terms of $(\rho^2)$-integers, that is, $(-\phi)$-integers, in the form
$$ \omega = \omega^{\prime} + \rho \, \omega^{\prime \prime}.$$

Suppose that $\Re(s)$ is sufficiently large.  Let $\Omega$ denote a set of  unique representations for the $\rho^2$-integers, and $\Omega_n \subset \Omega$ the representations of the  $n$-digit $\rho^2$-integers. Note that $\rho^2 = -\phi$. The particular set  $\Omega$ of unique $\rho^2$-integers, i.e., $(-\phi)$-integers, that we have in mind is  defined by the following relations:
\be
\left.
\ba{rcl}
\Omega_1 &=& \{ 0, 1 \} \\
\Omega_2 &=& \{ 0, 1, -\phi  \} \\
\Omega_n &=& ( - \phi)\, \Omega_{n-1} \cup ( 1+ (-\phi)^2 \, \Omega_{n-2} ) 
\ea
\right\}
\label{Omega-n decomp}
\ee
Then $\Omega = \bigcup_{n=1}^{\infty} \Omega_n = \lim_{n \rightarrow \infty} \Omega_n$.  Put $\Omega^* = \Omega - \{ 0 \}$.  Induction shows that $\#\Omega_n = F_{n+2}$.

Define the $\rho$-Weierstra{ss} $\wp$-function $\wp_{\rho}(z) $ and the related $\rho$-Hurwitz function $\wp_{\rho}^{(s)}(z)$ by\footnote{Notice that the sign convention departs from what is customary for the Hurwitz $\zeta$-function, in order to  position the poles of the meromorphic function at the $\rho$-integer points of the  associated tessellation.}
\bea
\wp_{\rho}(z) &=& \frac{1}{z^2} + \sum_{\omega \in \Omega^*}   \left\{   \frac{1}{(z - \omega)^2} - \frac{1}{(- \omega)^2}  \right\}  \\
\wp^{(s)}_{\rho}(z) &=&  \sum_{\omega \in \Omega }  \frac{1}{(z - \omega)^s}  , \quad \Re(s) >2 .
\eea
$\wp^{(n)}_{\rho}(z)$ is proportional to the $n$-th derivative of $\wp_{\rho}(z)$.

In terms of $\Omega$, the  $\rho$-Weierstra{\ss} function is
\be
 \wp_{\rho}(z) = \frac{1}{z^2} + 
 	\sum_{\omega^{\prime}, \omega^{\prime \prime}  \in \Omega}
	 \left( 
	 \frac{1}{z-(\omega^{\prime} + \rho \, \omega^{\prime \prime})^2} -  
	 \frac{1}{(\omega^{\prime} + \rho \, \omega^{\prime \prime})^2} 
	 \right) ,
 \ee
 with a corresponding expression for  $\wp_{\rho}^{(s)}(z)$. We shall show that these sums converge in section \ref{sec: rho in C}.
 
\begin{figure}[t]
\begin{center}
\caption{Radix $\rho = i \sqrt{\phi} $. The modulus of the $\rho$-Weierstra{\ss} $\wp$-function. 6 digit series expansion of the function. The domain is the  remainder set  $R$. The decomposition of  $R$  into four essentially disjoint tiles $R_1 \cup R_2 \cup R_3 \cup R_4$ is evident (Cp. fig.\ref{HLR_2D-rtPhi_R_decomp}, page \pageref{HLR_2D-rtPhi_R_decomp}).
}
\label{P rho=i rt phi dig=6} 
\includegraphics[width=3in]{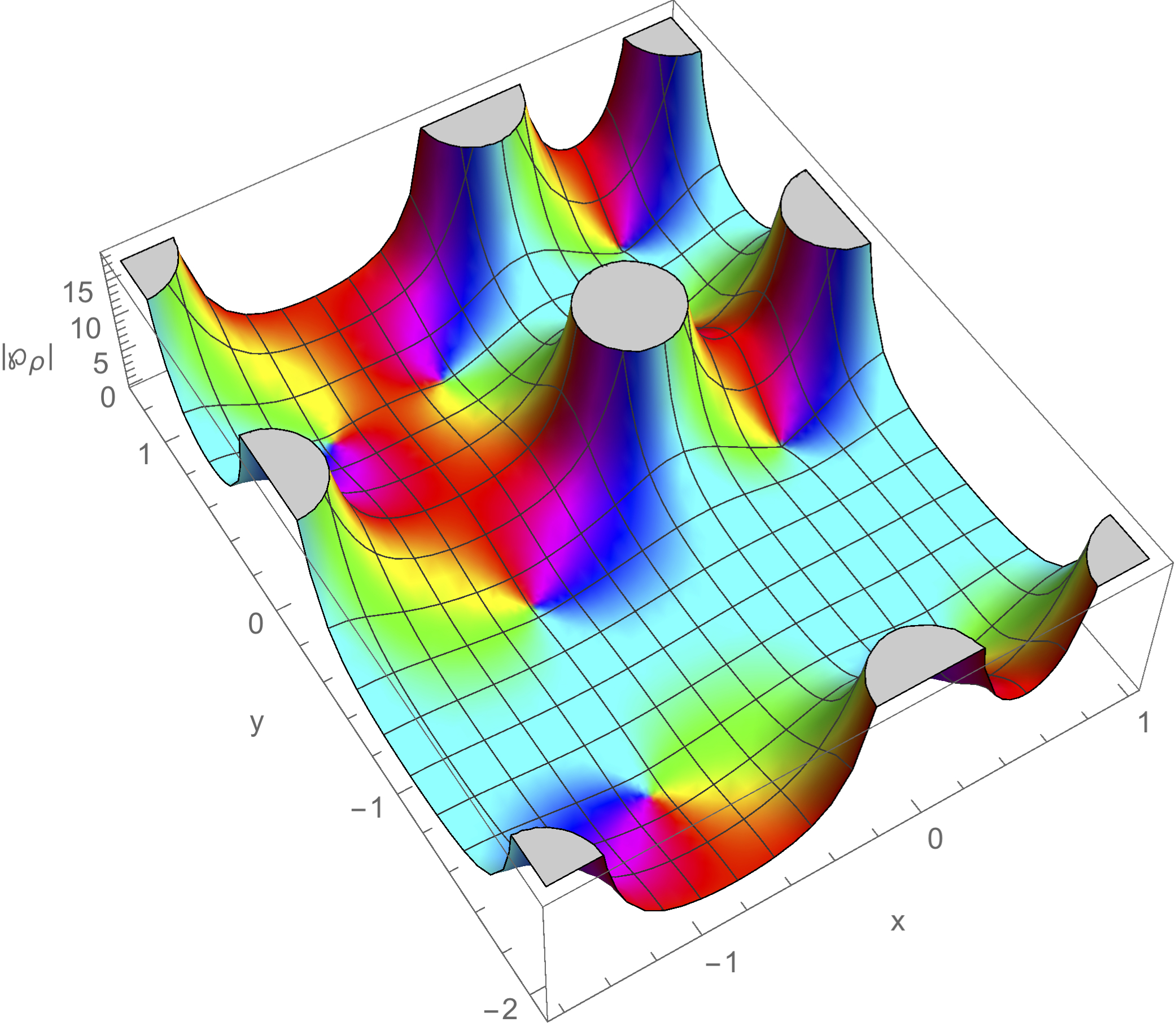}
\end{center}
\end{figure}

\begin{figure}[b]
\begin{center}
\caption{Radix $\rho = i \sqrt{\phi} $. The modulus of $ d \wp_{\rho}(z)/ dz $. 6 digit series expansion of the function. The domain is the  remainder set  $R$. The decomposition of  $R$  into four essentially disjoint tiles $R_1 \cup R_2 \cup R_3 \cup R_4$ is evident (Cp. fig.\ref{HLR_2D-rtPhi_R_decomp}, page \pageref{HLR_2D-rtPhi_R_decomp}).
}
\label{P rho=i rt phi s=3 dig=6} 
\includegraphics[width=2.5in]{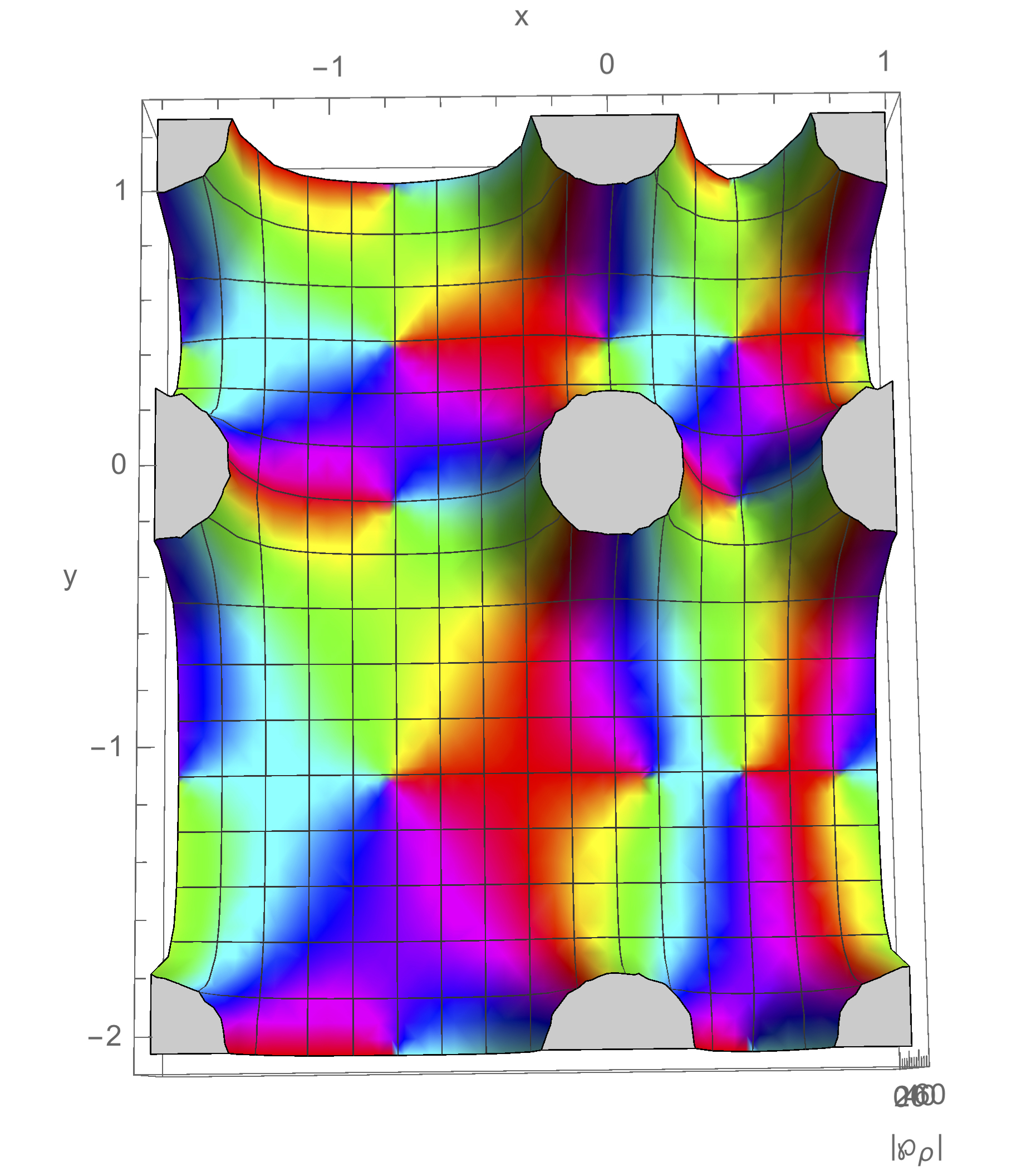}
\end{center}
\end{figure}

Figures \ref{P rho=i rt phi dig=6} and \ref{P rho=i rt phi s=3 dig=6} display the modulus of $ \wp_{\rho}(z)$  and $ \wp_{\rho}^{\prime}( z)$ respectively.

\aster
\clearpage

In constructing the $\rho$-Weierstra{\ss} function and deriving its $\rho$-duplication formula we shall make use of helping functions $\wp_{\rho,n_1,n_2}(z)$ and $\wp^{(s)}_{\rho,n_1,n_2}(z)$  that are finite sums whose limits will be $\wp^{(s)}_{\rho}(z) $. Introduce the finite sum
\be
\wp_{\rho,n_1,n_2}^{(s)}(z) :=  \sum_{\substack{ \omega^{\prime}_{n_1}  \in \Omega_{n_1} \\ \omega^{\prime \prime}_{n_2} \in \Omega_{n_2} }} \frac{1}{ \left( z - ( \omega_{n_1}^{\prime}  + \rho \,  \omega_{n_2}^{\prime \prime} ) \right)^s}.
\ee
and the corresponding finite sum for $\wp_{\rho,n_1,n_2}(z)$.
If the limit function exists, we have
\be
 \wp_{\rho}^{(s)}(z) = \lim_{\substack{n_1 \rightarrow \infty \\ n_2 \rightarrow \infty}} \wp_{\rho,n_1,n_2}^{(s)}(z),
 \ee
 with the corresponding limit  formula for $\wp_{\rho,n_1,n_2}(z)$.

\vspace{12pt}

\begin{it}
The following identity for finite sums is the $\rho$-duplication formula for $\wp_{\rho,n,n}^{(s)}(z)$: 
\bea
\rho^s \, \wp_{\rho,n,n}^{(s)}(\rho z) 
	&=& \frac{1}{\rho^s}  \wp_{\rho,n-1,n-1}^{(s)}\left( \frac{z}{\rho} \right)   		+  \frac{1}{\rho^{2s}}  \wp_{\rho,n-1,n-2}^{(s)}\left( \frac{z}			{\rho^2} - \frac{1}{\rho^3} \right)   +\nonumber \\
	&& \qquad  \frac{1}{\rho^{3s}}  \wp_{\rho,n-2,n-2}^{(s)}\left( \frac{z}		{\rho^3} - \frac{1+\rho}{\rho^4} \right)   +   \frac{1}{\rho^{3s}}  		\wp_{\rho,n-2,n-2}^{(s)}\left( \frac{z}{\rho^3} - \frac{1}{\rho^3} 		\right)   + \nonumber \\
	&& \qquad \qquad  \frac{1}{\rho^{4s}}  \wp_{\rho,n-2,n-3}^{(s)}			\left( \frac{z}{\rho^4} - \frac{1+\rho}{\rho^4} \right)    .
\label{finite sum p-fncn Hurwitz s dup}
\eea
\end{it}

\pf The proof consist of a verification that the various terms on the left side of eq(\ref{finite sum p-fncn Hurwitz s dup}) can be re-written as the terms on the right side. This is accomplished by using eq(\ref{Omega-n decomp}). The general summand on the left is
$$
\left( 
z - \frac{
	\omega^{\prime}_n + \rho \, \omega^{\prime \prime}_n
	}{\rho} 
\right)^{-s} .
$$
From eq(\ref{Omega-n decomp}) we see that  one way  $\omega_n$ can be expressed is $\rho^2 \, \omega_{n-1}$, Apply this to $\omega^{\prime}_n,  \omega^{\prime \prime}_n$ to find
$$
 \frac{
	\omega^{\prime}_n + \rho \, \omega^{\prime \prime}_n
	}{\rho} 
  \rightarrow 
\rho \, 
\left( 
	\omega^{\prime}_{n-1} + \rho \, \omega^{\prime \prime}_{n-1}
\right) ,
$$ 
and move the factor $\rho$ outside the sum to obtain  $ \rho^{-s} \wp_{\rho,n-1,n-1}^{(s)}\left( z/\rho \right) $, which is the first term in eq(\ref{finite sum p-fncn Hurwitz s dup}). The other terms are handled similarly. $\Box$ \\

\begin{it}
If the series converges, the $\rho$-duplication formula for $\wp_{\rho}^{(s)}(z)$ is
\bea
\rho^s \, \wp_{\rho}^{(s)}(\rho z) 
	&=& \frac{1}{\rho^s}  \wp_{\rho}^{(s)}\left( \frac{z}{\rho} \right)   			+  \frac{1}{\rho^{2s}}  \wp_{\rho}^{(s)}\left( \frac{z}{\rho^2} - 		\frac{1}{\rho^3} \right)   + \frac{1}{\rho^{3s}}  \wp_{\rho}^{(s)}		\left( \frac{z}{\rho^3} - \frac{1+\rho}{\rho^4} \right)   +  				\nonumber \\
	&&  \frac{1}{\rho^{3s}}  \wp_{\rho}^{(s)}\left( \frac{z}{\rho^3} - 
		\frac{1}{\rho^3} \right)   +  \frac{1}{\rho^{4s}}  \wp_{\rho}^{(s)}		\left( \frac{z}{\rho^4} - \frac{1+\rho}{\rho^4} \right)    .
\eea
\label{rho-wp s-deriv dup}
\end{it}

\pf Take the limit as $n \rightarrow \infty$ in the previous theorem. It will converge, as we shall later see, for sufficiently large $\Re(s)$. $\Box$\\

The derivation of eq(\ref{rho-wp s-deriv dup}) did not depend on $s$ being an integer; it sufficed that $\Re(s)$ be large enough for the underlying series to converge. Therefore, the $\rho$-duplication formula is valid for the $\rho$-Hurwitz $\zeta$-function for arbitrary $s \in \C$ with large enough real part. We shall use this fact below. \\

Regarding  $\wp_{\rho}(z)$ itself, application of the same method yields:\\

\begin{it}
 The  $\rho$-duplication formula  for $\wp_{\rho}(z)$ is
\bea
\rho^2\, \wp_{\rho} \left( \rho z \right) 
	&=& 
	 \frac{1}{\rho^2} \, \wp_{\rho} \left( \frac{z}{\rho} \right) +
	  \frac{1}{\rho^4} \, \wp_{\rho} \left(\frac{z}{\rho^2} -\frac{1}{\rho^3} 	\right) + \nonumber \\
	& & \hspace{-40pt}  \frac{1}{\rho^6} \, \wp_{\rho} \left(\frac{z}{\rho^3} 	-\frac{1+\rho}{\rho^4} \right)  +  \frac{1}{\rho^6} \, \wp_{\rho} 		\left(\frac{z}{\rho^3} -\frac{1}{\rho^3} \right) +   \frac{1}{\rho^8} \, 		\wp_{\rho} \left(\frac{z}{\rho^4} -\frac{1+\rho}{\rho^4} \right)   + C 
\label{rho-wp dup}
\eea
where
\bea
C &=&   \left(
		 \frac{1}{\rho^4} \, \wp_{\rho} \left( -\frac{1}{\rho^3} \right) +   		\frac{1}{\rho^6} \, \wp_{\rho} \left(  -\frac{1+\rho}{\rho^4} \right)  		+  \frac{1}{\rho^6} \, \wp_{\rho} \left(  -\frac{1}{\rho^3} \right) +   		\frac{1}{\rho^8} \, \wp_{\rho} \left(  -\frac{1+\rho}{\rho^4} \right) 
	\right) .
\eea

\end{it}



\section{Convergence properties}
\label{sec: convergence}

It is time to take up the question of convergence of the  infinite series we have been freely manipulating. In the framework of this  paper, we shall not be concerned with establishing the  tightest conditions for convergence; rather, we merely want to show that the series converge for some potentially significant values of the variables so that the functions we have introduced, and the identities they satisfy, make sense.

The convergence problems concern the $\rho$-Riemann $\zeta$-function;  the $\rho$-Hurwitz $\zeta$-function; and the $\rho$-Weierstra{\ss} $\wp$-function. In the classical case one first proves convergence of the series that defines the Riemann $\zeta$-function, reduces the convergence of the Hurwitz $\zeta$-function to it; and applies the same approach to the Weierstra{\ss} $\wp$-function, which  can be thought of as a 2-dimensional version of the Hurwitz function.  

The new series are  sums over a set of $\rho$-integers. The convergence proofs follow the general pattern of the corresponding classical case; cp., for instance, the famous textbook of Whittaker and Watson \cite{Whittaker+Watson}. The principal difference is that the summands are collected according to the leading term in the positional representation of the $\rho$-integer summation variable.

As above, let $1< | \rho | \leq 2$ be a radix and $\Delta=\{ 0, 1 \}$,  and let $\Delta[ \rho ]$ denote the set of $\rho$-integers. In general, the notation $ \sum_{S^*}$  stands for the sum over the index set $S^*:=S-\{0\}$.  The convergence of the series will be subdivided into cases, according to the nature of $\rho$: $\rho \in \R^+$; $\rho \in \R^-$; and $\rho \in \C - \R^+$.  Within each case we shall also have to consider whether each $\rho$-integer has a unique positional representation.

The idea behind the proof of convergence is to collect terms in the sum that have the same highest power of $\rho$. Let $\Delta[\rho]_n$ denote the set of $n$-digit $\rho$-integers. The set of $\rho$-integers whose highest power of the radix is $\rho^n$ is $\nabla_n := \Delta[\rho]_n -\Delta[\rho]_{n-1}$.

Begin by considering the $\rho$-Riemann $\zeta$-function; cp. eq(\ref{rho-Riemann def}).
$$
 \zeta_{\rho}(s) := \sum_{\omega \in \Delta [ \rho ]^*} \frac{1}{\omega^s} .
$$

\subsection{$\rho \in \R^+$}
\label{sec: R+ conv}

This is the simplest case because a $\rho$-integer is a sum of non-negative monomials. 

Let $\omega_n =  \rho^{n-1} + \sum_{k=0}^{n-2} \delta_k \rho^k  \in \nabla_n$ be a $\rho$-integer whose highest power of the radix is $n$.   Then
$$ \zeta_{\rho}(s) = \sum_{n=1}^{\infty} \sum_{\nabla_n} \frac{1}{\omega_{n}^s} $$
and
$$ 
\left| \zeta_{\rho}(s) \right| 
	\leq  \sum_{n=1}^{\infty} \sum_{\nabla_n}\frac{1}{  \left|  \omega_{n}^{s}  \right| }
		\leq   \sum_{n=1}^{\infty} \sum_{\nabla_n}  \frac{1}{\omega_{n}^{\Re(s)} } .
		$$
Evidently
$$ \rho^{n-1} \leq \omega_n \leq \frac{\rho^n -1}{\rho -1}.$$
Since the $\delta_k$ are 0 or 1, there are at most $2^{n-1}$ distinct $n$-digit $\rho$-integers, whence
$$
2^{n-1}  \frac{\rho -1}{\rho^n -1} \leq 2^{n-1} \frac{1}{\omega_n} \leq \frac{2^{n-1}}{\rho^{n-1}} ,
$$
so
\be
(\rho -1)^{\Re(s)} \sum_{n=1}^{\infty} \frac{2^{n-1}}{\rho^{n-1} -1} 
	\leq \left| \zeta_{\rho}(s) \right| 
	\leq \sum_{n=1}^{\infty} \left( \frac{2}{\rho^s}\right)^{n-1} 
	= \frac{1}{1 - 2 \rho^{-s}} ,
\ee
so we have proved that if $1 < \rho \leq 2$ and $\Delta = \{ 0,1 \}$,  then  $ \zeta_{\rho}(s) $ converges absolutely for $\Re(s)>  \log 2/ \log \rho$.

Note that $\Re (s) >  1$ for all radices. Recall that the representation of $\rho$-integers need not be unique. Nevertheless, the lemma applies, for instance, to $\rho = \phi$ where positional representation of $\rho$-integers is not unique because the radix is positive. In this case, we have convergence for $\Re(s) > \log2 /\log \phi \simeq =1.44 \dots .$

\aster

The radix $\rho = \phi$ generates more than one representation for some $\rho$-integers. As we  show in section \ref{sec:radix phi}, the cardinality of $\Omega_n$, i.e. the number of  pairwise unequal representations for $n$-digit $\phi$-integers is  $F_{n+2}$ where $F_n$ is the Fibonacci number. The convergence criterion for 
$$\zeta_{\rho}(s) = \sum_{\Omega^*} \omega^{-s}$$
 is derived as above, but now the cardinality of  $ \nabla_n$ is $F_{n}$ rather than $2^n$ so $| \zeta_{\rho}(s) | \leq \sum_{n=1}^{\infty}  \frac{F_{n}}{\phi^{(n-1) \Re(s)} }$. Since $F_{n}$ is the closest integer to $\phi^{n}/\sqrt{5}$, the series converges when
 $$ 
\sum_{n \geq 1} \frac{F_{n}}{(\phi^{\Re(s)})^{n-1} }  \sim \sum_{n \geq 1}     \left( \frac{1}{\phi^{\Re(s)-1} } \right)^{n-1}
 $$ 
converges,  that is, when $\Re(s) > 1$. This improves on the general condition $\Re(s) > \log 2/\log \rho$ given above.

\aster

Convergence of the $\rho$-Hurwitz is dealt with in the usual way, by noting that 
$$
\left| \sum_{\Omega} \frac{1}{(z + \omega)^s} \right| \leq 
	\sum_{\Omega} \frac{1}{\left| z + \omega \right|^{\Re(s)}},
$$
and that for $z$ in a compact set, $\left| z + \omega \right|$, apart from finitely many terms, the terms of the series are dominated by the corresponding terms of $\zeta_{\rho}(s)$, for which convergence has already been shown. 

\subsection{$\rho \in \R^-$}
\label{sec: rho in R}

Now suppose that $-2 \leq \rho <-1$. The $\rho$-integers have negative as well as positive values that are unbounded as the number of digits is increased. Again consider an $n$-digit $\rho$-integer whose leading coefficient is 1.

$$ 
\omega_n = \sum_{k=0} \delta_{2k} \rho^{2k} + \rho \,  \sum_{k=0} \delta_{2k+1} \rho^{2k},
$$
where the sums are non negative. $\omega_n$ is positive, resp., negative, if  $n-1$ is even, resp., odd.  If $n-1 = 2l$, then
$$
\rho^{2l} + \rho \left( \frac{ \rho^{2l} -1}{\rho^2 -1} \right) \leq  \omega_n   \leq \frac{\rho^{2l+2} -1}{\rho^2-1}.
$$
A $\rho$-integer that is closest to the origin is $\rho+1$. If $\rho \neq -2$, it is unique. Let $\omega_n$ denote, as above, a $\rho$-integer that is expressed by exactly $n$-digits, so that the coefficient of $\rho^{n-1}$ is 1. 

For $\rho < -1$ we cannot hope to bound $| \omega_n |$ away from 0 by a bound that increases with $n$ in general. For example, if $\rho = -\left( \frac{1+\sqrt{5}}{2} \right)$, then
 $$ \rho^{2n}+ \sum_{k=1}^n \rho^{2k-1}=1, $$
 so the set of $\rho$-integers with an odd number of digits always contains 1. This implies that sums such as $\sum_{\Delta[\rho]^*} \omega^{-s}$ will have an unlimited number of terms equal to $1^{-s}$ and cannot converge.  This shows that a radix for which the representation of a $\rho$-integer is not unique require special, and perhaps individual, consideration to eliminate duplications.
 
We restrict consideration to $\rho = -\phi$ and to summations over the set $\Omega^0$ of even $\rho$-integers defined by eq(\ref{-phi V_n recursion}). Recall that $\#\Omega^0_n = F_{n+1}$.  Each $\omega \in \Omega^0_n$ has the form
$$ \omega = \sum_{k=0}^{n-1} \delta_k \rho^k . $$

Denote the subset of elements of $\Omega_n^0$ whose leading coefficient $\delta_{n-1}=1$ by $\nabla^0_n$.  Thus $\nabla^0_n=\Omega^0_n-\Omega^0_{n-1}$. Evidently $\#\nabla^0_n = F_{n-1}$. The non-zero elements of $\nabla^0_n$ are negative, resp., positive, if $n$ is even, resp., odd. 

We seek a lower bound for $\omega_n \in \nabla^0_n $, and find $F_{n} <  | \omega_n |$. Thus
\bean
\left|   \sum_{(\Omega^0_n)^* } \frac{1}{ \omega^s} \right|  & < & 
 	 \sum_{(\Omega^0_n)^*} \frac{1}{ \left| \omega^s  \right|} \\
	 	& =  &  \sum_{n=1}^{\infty} \sum_{(\nabla^0_n)^*} \frac{1}			{ \left| \omega^s  \right|} \\
		& =  &  \sum_{n=1}^{\infty} \sum_{(\nabla^0_n)^*} \frac{1}			{ (F_{n-1})^s} \\
		&=&  \sum_{n=1}^{\infty} \frac{1}{ (F_{n-1})^{s-1}} \\
		& \sim &  \sum_{n=1}^{\infty} \frac{1}							{ \left( \phi^{s-1}\right)^n }\\
\eean
whence the series converges for $\Re(s) > 1 + \frac{1}{\log \phi} \sim 3 $.

\subsection{$\rho \in \C -  \R$ with $\Delta =   \{ 0,1 \}$}
\label{sec: rho in C}

We confine ourselves to the complex radices that are the classical  complex multiplications for elliptic functions, and to $\rho = i \sqrt{\phi}$.  

In the classical cases,  where $\rho$ is an imaginary quadratic integer, an arbitrary $\rho$ integer has the form $m + n \rho$, where $m,n \in \Z$.\footnote{The $\rho$-Riemann $\zeta$-function $\zeta_{\rho}(s) = \sum_{\Omega^*} \omega^{-s}$  can be interpreted as a special value of an Eisenstein series when $s \in 2\Z^+$.}  For these cases, the problem of convergence proceeds as usual: a sequence of artfully arranged similar nested  parallelograms covers the non-zero $\rho$-integers; a counting argument bounds the number terms on each parallelogram; and an estimate for a lower bound is found that is proportional to the index of the parallelogram. 

In the general 2-dimensional case we have neither the 2-dimensional basis nor the parallelograms. Instead, we partition $\Omega^*$ into sets $\nabla_n$ as we have on other occasions. In order to do this we must make explicit the structure of a complete set of unique representatives $\Omega$ for the $\rho$-integers. 

For $\rho =i \sqrt{\phi}$ the geometry is much simpler, and this is why we use it as an example. The $\rho$-integers are arranged on a sequence of concentric rectangles centered on the origin. They are not distributed in the regular manner found in the classical case, but they are nevertheless easy to count.

Number successive similar rectangles starting with 1, and denote the $k$-th rectangle by ${\rm rect}_k$. The number  of $\rho$-integers on  ${\rm rect}$ is $ 8 k$ (cp. fig.\ref{dim2 silver(0101) R1 10 digit}) and the total number of $\rho$-integers, including $0$,  in the closure of the $k$-th rectangle is $1+ 8 \sum_{n=1}^{k} n = (2k+1)^2$. These rectangles are successive inflations of a remainder set for radix $\rho$, which provides the foundation for a detailed proof.

The minimum distance from the origin to ${\rm rect}_k$ is achieved along the positive real axis. The differences between successive real $\rho$-integers is either $1$ or $\phi$, whence the distance of points on ${\rm rect}_k$ from 0 is greater than or equal to $k$. This will suffice to show that the various series converge for interesting values of the variable.

Consider. for instance, the $\rho$-Riemann $\zeta$-function,
$$ \zeta_{\rho}(s) = \sum_{\Omega^*}  \frac{1}{\omega^s} .$$
We find
\bean
\left|  \zeta_{\rho}(s) \right| & = &
	\left|  \sum_{\Omega^*} \frac{1}{\omega^s } \right| 
	\leq  \sum_{k=1}^{\infty} \sum_{{\rm rect}_k} \frac{1}{\left|  			\omega_k\right|^{\Re(s)} }  
	\leq   \sum_{k=1}^{\infty}   \frac{8k}{k^{\Re(s)} } \\
	& \sim &  \sum_{k=1}^{\infty}   \frac{1}{ k^{\Re(s)-1}  },
\eean
which converges for $\Re(s)>2$.  This simple result is enough to conclude that the series for the derivatives of the $\rho$-Weierstra{\ss} $\wp_{\rho}$-function converge, because they are just the numbers $\zeta_{\rho}(k+2)$, $k \in \Z^+$. 

Regarding $\wp_{\rho}(z)$ itself, we follow the standard proof. Suppose that $z$ is an element of a compact subset of $\C - \Omega$.  The absolute value of the general summand satisfies
\bean
\left| \frac{1}{(z - \omega)^2} - \frac{1}{ \omega^2} \right| & = &
	|z| \, \frac{| z-2 \omega | }{| z-\omega |^2 | \omega |^2} \\
	& \sim & \frac{1}{|\omega|^3},
\eean
whence the series for $\wp_{\rho}(z)$ is  uniformly convergent on compact subsets of $\C - \Omega$.

\subsection{Analytic continuations}
\label{sec: analytic continuations}

Equation (\ref{zeta_rho-Riemann Laurent series}) and its special case eq(\ref{analytic continuation: classical}) for $\rho=2$, and  eq(\ref{phi zeta Laurent series}) for $\rho = \phi$, provide an analytic continuation  of the respective $\zeta_{\rho}(s)$ to $\C$. Comparison shows that the two expressions have a similar form, apart from a common factor that depends only on the radix and $s$, the power of $\rho$ that appears in the denominator of each summand, and the intrinsic difference between the two $\rho$-Riemann $\zeta$-functions. Here we shall investigate convergence of the shared series
\be
1 + \sum_{k \geq 1} (-1)^k \frac{\Gamma(s+k) }{ \Gamma(s) \Gamma(k+1) } \frac{ \zeta_{\rho}(s+k) }{ \rho^{c(s+k)}} ,
\ee
where $c=1$ or $c=2$.

We have shown that $\zeta_{\rho}(s)$ converges for $\Re(s)$ sufficiently large, and it is evident that  for any $\epsilon > 0$, there is an $s_{\epsilon}$ such that  $\left| \zeta_{\rho}(s) -1  \right| < \epsilon$; hence, as far as convergence is concerned, we can ignore the factor $\zeta_{\rho}(s+k)$.  

Temporarily write $s$ for $\Re(s)$ to simplify notation. Application of Stirling's formula shows that
$$ 
\frac{\Gamma(s+k) }{ \Gamma(s) \Gamma(k+1) } \frac{ 1 }{ | \rho |^{c(s+k)}} \sim \frac{e^{s-1} }{s^{s+1}} \frac{(1+s/k)^{s+1}}{(1+1/k)^2} \frac{k^{s-1} }{ | \rho |^{c(s+k)}} \sim \frac{k^{s-1} }{ | \rho |^{c(s+k)}}.
$$
The terms go to 0 with increasing $k$. The comparison test shows that the ratio of successive terms is ultimately, i.e., for sufficiently large $\Re(s)$, $1/| \rho|^c \leq 1/|\rho| < 1$, so the series converges to the right of the pole give by the factor preceding the series. This factor was $ \phi^{s}/(\phi^{s} -1)$ if the representative polynomials $\Delta [\rho]$ are unique,  and $\phi^{2s}/(\phi^{2s} - \phi^{s} -1)$ for $\rho =\phi$.

\section{Remarks}
\label{sec:remarks}

\subsection{About the color coding of the graph of a complex valued function}
\label{sec: color coding}

The graph of the modulus $f(z)$ of a complex valued function is colored to show the value of  the argument of $f(z)$. The colors are drawn from the circle of hues so that positive real numbers (argument $=0$)  correspond to red. As the argument increases from 0 to $2 \pi $, the hues range from red through yellow, green, teal (negative reals) blue, and magenta.

\subsection{About the numerical calculation of identities}
\label{sec: numerics}

I have checked many of the new -- especially  the duplication formulae --  numerically. The calculations are numerically unstable because they require cancellation of values in  neighborhoods of the many poles of the functions. A few remarks about how this can be done may be helpful.

The problem can be illustrated by examining expressions like $f(s,z) = (z+ \omega)^{-s}$  that appear in many of the formulae. The $\rho$-duplication formulae involve expressions such as 
$$
\rho^s \, f(s, \rho z), \quad  f(s, z), \quad \frac{1}{\rho} f(s, z/\rho + 1/\rho^2).
$$
which the theorem tells us cancel.

Because all numerical calculations are performed with a fixed maximum  number of digits, it is essential to insure that only terms with the same number of digits appear. Recall that $\Omega_n$ denotes the set of $n$-digit $\rho$-integers for some radix  $\rho$. Set $f_n(s,z) = (z+ \omega_n)^{-s}$ where $\omega_n \in \Omega_n$.  
Multiplication and division of $z$ by powers of the radix changes the number of digits. For instance, 
$$ \rho^s \, f_n(s,\rho z) =\frac{1}{(  z + \omega_n/\rho)^s} $$
so the numbers $\omega_n/\rho$ for which the low order digit of $\omega_n$ is zero  actually belong to $\Omega_{n-1}$ rather than $\Omega_n$. Identities such as those in section \ref{numerical exact identities}  insure that the terms in a calculation are effectively of the same order, thereby enabling cancellations to an arbitrarily high degree of accuracy.

\noindent 20150724 \\
\noindent Boston and Gloucester

\end{document}